\documentclass[12pt]{article}
\usepackage{amsmath}
\usepackage{amssymb}

\usepackage{graphicx}

\let\rho=\varrho

\newcommand{\qed}{\hfill$\Box$\par\medskip\par\relax}

\newcommand{\eps}{\varepsilon}
\newcommand{\Z}{{\mathbb Z}}
\newcommand{\N}{{\mathbb N}}
\newcommand{\V}{{\mathcal V}}
\newcommand{\K}{{\mathcal K}}
\newcommand{\hK}{{\hat{\mathcal K}}}
\newcommand{\M}{{\mathcal M}}

\newcommand{\HH}{{\mathcal H}}
\newcommand{\B}{{\mathcal B}}
\newcommand{\R}{{\mathbb R}}

\newcommand{\ZZ}{{\mathcal Z}}

\newcommand{\Q}{{\mathbb Q}}
\let\phi=\varphi
\newcommand{\A}{{\mathfrak A}}

\newcommand{\br}{{\mathsf b}{\mathsf r}}

\newcommand{\tV}{{\tilde{\mathcal V}}}

\newcommand{\eqlaw}{\stackrel{\text{\tiny law}}{=}}

\newcommand{\om}{\text{\boldmath ${\omega}$}}

\newcommand{\dist}{{\mathop{\rm dist}}}

\newcommand{\EA}{{\mathbf E}}
\newcommand{\PA}{{\mathbf P}}
\newcommand{\IE}{{\mathbb E}}
\newcommand{\IP}{{\mathbb P}}
\newcommand{\Po}{{\mathtt P}_{\omega}}
\newcommand{\Eo}{{\mathtt E}_{\omega}}
\newcommand{\PoA}[1]{{\mathtt P}_{\omega|#1}}
\newcommand{\EoA}[1]{{\mathtt E}_{\omega|#1}}
\newcommand{\p}{{\mathbf p}}

\newcommand{\8}{{\infty}}

\newcommand{\lf}{\lfloor}
\newcommand{\rf}{\rfloor}

\newcommand{\Do}{\Delta^{\!\omega}}

\newcommand{\I}{{\mathbf 1}}
\newcommand{\1}[1]{{\bf 1}{\{#1\}}}
\newcommand{\BRW}{{branching random walk}}

\allowdisplaybreaks

\newtheorem{theo}{Theorem}[section]

\newtheorem{lmm}[theo]{Lemma}
\newtheorem{df}[theo]{Definition}
\newtheorem{prop}[theo]{Proposition}

\newtheorem{rem}[theo]{Remark}
\newtheorem{example}[theo]{Example}

\title{Shape and local growth for multidimensional branching
      random walks in random environment}

\author{Francis~Comets\thanks{Partially
   supported by CNRS (UMR 7599 ``Probabilit{\'e}s et Mod{\`e}les
Al{\'e}atoires'') and by the 
``R\'eseau Math\'ematique France-Br\'esil''}$^{~,1}$ \and
 Serguei~Popov\thanks{Partially supported by CNPq (302981/02--0 and 471925/2006--3),
FAPESP (04/07276--2), USP/COFECUB, and by the ``Rede Mate\-m\'atica
Brasil-Fran\c{c}a"}$^{~,2}$}

\begin{document}

\maketitle

{\footnotesize
\noindent $^{~1}$Universit{\'e} Paris 7, UFR de Math{\'e}matiques,
case 7012, 2, place Jussieu, F--75251 Paris Cedex 05, France

\noindent e-mail: \texttt{comets@math.jussieu.fr}, \quad
\noindent url: \texttt{http:/\!/www.proba.jussieu.fr/$\sim$comets}

\noindent $^{~2}$Instituto de Matem{\'a}tica e Estat{\'\i}stica,
Universidade de S{\~a}o Paulo, rua do Mat{\~a}o 1010, CEP 05508--090,
S{\~a}o Paulo SP, Brasil

\noindent e-mail: \texttt{popov@ime.usp.br}, \quad
\noindent url: \texttt{http:/\!/www.ime.usp.br/$\sim$popov}

}

\begin{abstract}
We study \BRW s in random  environment on the $d$-dimensional square lattice,
$d \geq 1$.
In this model, the environment has finite range dependence, and 
the population size cannot decrease.
We prove limit theorems (laws of large numbers) for the set of lattice
sites which are visited up to a large time as well as for the local size
of the population. The limiting shape of this set is compact and convex,
and the local size is given by a concave growth exponent.
Also, we obtain the law of large numbers for
the logarithm of the total number of particles in the process.
\\[.3cm]{\bf Short Title:}  Branching random walks in random environment
\\[.3cm]{\bf Keywords:} shape theorem, 
subadditive ergodic theorem, transience,  growth exponent, population size
\\[.3cm]{\bf AMS 2000 subject classifications:} Primary 60K37;
secondary 60J80, 82D30
\end{abstract}

\section{Introduction and results}
\label{s_intro}
We start with an informal description of
%to describe informally 
the model we study in this paper.
Particles live in~$\Z^d$ and evolve in discrete time.
At each time, every particle is substituted by (possibly more than one)
offspring which are placed in neighboring sites, independently of the other particles.
The rules of offspring generation depend only on the location of the particle.
The collection of those rules (so-called the \emph{environment}) is itself 
random,
it is chosen randomly before starting the process, and then it is kept
fixed during
all the subsequent evolution of the particle system.

This model considered in this paper was introduced
in~\cite{CP}. The random environment here affects both branching
and transition mechanisms, and (as opposed to the models of~\cite{CMP,GM,MP1,MP2})
the \emph{immediate} descendants of a particle are not supposed to be independent.
%%    In~\cite{CP} this model was
%%     studied on the grounds of recurrence and transience, shape theorems,
%%    hitting times. 
In~\cite{CP} we proved a dichotomy for  recurrence/transience, depending
only on the support of the environmental law, and we gave sufficient conditions
for  recurrence and transience. In the recurrent case, the tails of the 
hitting times are studied and a shape theorem is obtained in a preliminary
form.
The recurrence/transience classification was later completed in~\cite{M1,M2}. 
We refer also to \cite{BCGdH,Eng,HMP,Volk} for other work related to
branching random walks in random environment.

 Now, let us look at the subject of this paper from a different angle.
For  $x\in \Z^d$, let  $p(x,\cdot)$ be the transition probabilities
 from $x$ to its nearest neighbors, and %%``branching'' number
$r(x) \geq 1$. Consistently with the notation introduced later in this 
section,
we denote by $\om=(\omega_x; x \in \Z^d)$ the collection of coefficients 
$\omega_x=(p(x,x+\cdot), r(x))$ (as explained below,
$r(x)$ stands for the mean offspring in~$x$), and by~$\Do$ the corresponding discrete Markov
%% Laplacian NON REVERSIBLE
operator, 
\[
\Do f(x)= \sum_{y \sim x} p(x,y)[f(y)-f(x)].
\]
In this paper we will study for $z \in \Z^d$ 
the solution $u_n(x)=u_n^z(x)$ of the equation
\begin{equation}  
\label{eq:discrete_anderson}
\left\{
\begin{array}{ll}
u_{n+1}-u_n= r \Do u_n + (r-1) u_n, &  x \in \Z^d, \;
n =0,1,\ldots,\\
u_0(x)=\1{x=z}\vphantom{\sum^N}. \\
\end{array} 
\right.
\end{equation}
It is easily checked, for instance 
by the discrete
Feynman-Kac formula, that the solution $u_n^z$ 
is given by the expectation
$u_n^z(x)=\Eo \eta_{n}^x(z)$ of the number $\eta_{n}^x(z)$ of particles 
in~$z$ at time~$n$ in a discrete-time \BRW~starting from a single 
individual located at site $x$ at time 0. The evolution rule of this \BRW{} 
is that particles at~$x$ branch with 
an average of~$r(x)$ children which then move independently to a neighboring
site randomly chosen from $p(x,\cdot)$. We will be interested in the case where
the coefficients~$\om$ are given by a
stationary and finitely dependent
%independent and identically distributed 
random field. The model has other possible 
formulations.
In the case of continuous time, the above equation becomes the parabolic 
partial differential equation
\[
\frac{\partial }{\partial t} u_t(x) = \Do u_t(x) + V^\omega(x) u_t(x)
\]
with $V^\omega(x)$ the branching rate. In our case the mean number of
offspring is greater than or equal to~$1$ and
bounded from above, or, equivalently, $V^\omega(x)$ is nonnegative and bounded.

In the case where $p(x,\cdot)$ are the simple random walk transition 
probabilities $(1/2d,\ldots,1/2d)$, so that $\Delta$ is the standard
Laplace operator, this equation is known as the parabolic Anderson problem, 
%%or the real random Schr\"odinger equation, 
and has also 
continuous-space versions, see~\cite{Sz} and~\cite{Sz_drift}. 
These models have motivated a huge scientific activity, with particular 
interest on localization and intermittency (e.g.~\cite{GKM, HKM})
and survival analysis~\cite{Antal},
leading to fine pictures in the different cases of bounded or unbounded~$V$'s.
We stress that the  Markov operator~$\Do$ is random in the 
present paper, a case that seems not to have been studied so far.
Moreover, $\Do$ is non symmetric, this makes
the model non-reversible with a possibly large drift, and
prevents the use of a spectral theory as in the above references.

%\medskip

The present paper represents a step towards the understanding of the 
equation~(\ref{eq:discrete_anderson}).
We will prove that the solution $u_n^z$ behaves like
$u_n^z(x) = \exp \{n \beta((z-x)/n)+o(n)\}$ as $n \to +\8$
(see Remark~\ref{rem_complete_shape} where we comment on the function~$\beta$).
We will even get into the finer study of the population 
 size $\eta_n^x(z)$ itself, and we 
prove the same asymptotics under the assumption that particles cannot die. 
Hence, the growth of the population is rather smooth at the coarse scale
$z-x={\mathcal O}(n)$. The occurrence of shape theorems and local growth rate
goes back at least to \cite{biggins} and \cite{laredo-rouault} for branching
random walks, and to \cite{greven-hollander} in random environment. 
In fact, our model is slightly more general than described above.
As mentioned in the beginning of this section, 
the branching and the displacement may be dependent. 
Also, the environment that we consider is not necessarily independent,
but we allow for a finite range dependence. We believe that the results
extend to more general dependence, but considering this case would require
an important additional technical work. 

\medskip

{\bf The model.} 
We now describe the model, keeping the notations of~\cite{CP} 
whenever possible.
Let $\Z_+=\{0,1,2,\ldots\}$ and $e_i$-s be the
coordinate vectors of $\Z^d$. We use the notation
$\|x\|=|x^{(1)}|+\cdots+|x^{(d)}|$
for $x=(x^{(1)},\ldots,x^{(d)})\in\R^d$ (or $x\in\Z^d$).
Define the distance between two sets $A,B\subset\R^d$ 
(or $A,B\subset\Z^d$) by 
\[ \dist(A,B) = \inf_{a\in A, b\in B} \|a-b\|.\]
Fix a finite set $\A\subset \Z^d$ such
that $\pm e_i \in \A$ for all $i=1,\ldots,d$. 
Define
\[
 \V = \Big\{v=(v_x, x\in \A) : v_x \in \Z_+, \sum_{x\in \A}v_x \geq 1 \Big\},
\]
and for $v\in\V$ put $|v|=\sum_{x\in \A}v_x$; 
note that $|v|\geq 1$ for all $v\in\V$.
Furthermore, let~$\M$ be the set of all probability
measures~$\omega$ on $\V$:
\[
 \M = \Big\{\omega = (\omega(v), v\in\V) : \omega(v)\geq 0 
\mbox{ for all }v\in\V,
                     \sum_{v\in\V}\omega(v)=1\Big\}.
\]
Then, suppose that $\om:=(\omega_x\in\M,x\in\Z^d)$ is a stationary ergodic
random field,
and denote by $\IP, \IE$ the probability and expectation with respect to~$\om$.
Throughout this paper we suppose that this field is also \emph{finitely dependent},
that is, the following condition holds: 

\medskip
\noindent
{\bf Condition I.} 
There exists a positive number $\rho$ such that 
for any two sets $A,B\subset\Z^d$ with $\dist(A,B)\geq \rho$
the following holds: the sigma-algebra generated by $(\omega_x, x\in A)$
is independent under $\IP$ from
the sigma-algebra generated by $(\omega_x, x\in B)$.
Note that $\rho=1$ corresponds to the case of independent identically 
distributed environment.

\medskip

The collection $\om = (\omega_x, x\in\Z^d)$ is called {\it the environment}.
Given the environment~$\om$, the evolution of the process is described in the
following way: start with one particle at some fixed
site of~$\Z^d$. At each integer time
the particles branch independently using the following mechanism:
for a particle at site $x\in\Z^d$, a random element $v=(v_y, y\in \A)$
is chosen with probability $\omega_x(v)$, and then the particle
is substituted by $v_y$ particles in $x+y$ for all $y\in \A$.
Note that the population never gets extinct, since every individual 
has at least one direct descendant.

For the process starting from one particle at~$x$, let us denote by
$\eta_n^{x}(y)$  the number of particles in~$y$ at time~$n$.
Define the random variable
\[
 \ZZ_n^x = \sum_{y\in\Z^d} \eta_n^x(y),
\]
i.e., $\ZZ_n^x$ is the total number of particles at moment~$n$ for the
process starting from~$x$.

We denote by $\Po^x, \Eo^x$  the (so-called ``quenched")
probability and expectation for the process starting from~$x$ in the fixed
environment~$\om$. We use the notation
$\PA^x[\,\cdot\,] = \IE\,\Po^x[\,\cdot\,]$
for the annealed law of the branching random walk in random environment, and
$\EA^x$ for the corresponding expectation. Also, sometimes we use the symbols
$\Po,\Eo,\PA,\EA$ without the corresponding superscripts when it can create
no confusion (e.g.\ when the starting point of the process is
indicated elsewhere).

Throughout this paper, and often without recalling it explicitly,
we suppose that the two conditions below are fulfilled:

\medskip
\noindent
{\bf Condition B.}
\[
\IP[\mbox{there exists } v\in\V \mbox{ such that }
                         \omega_0(v)>0 \mbox{ and } |v|\geq 2]>0.
\]

\medskip
\noindent
{\bf Condition UE.} For some $\eps_0>0$,
\[
 \IP\Big[\sum_{v:v_e\geq 1} \omega_0(v) \geq \eps_0 \mbox{ for any }
 e\in \{\pm e_1,\ldots,\pm e_d\} \Big] = 1.
\]
Condition B ensures that the model cannot be reduced to random walk without
 branching, and Condition~UE is a natural uniform
ellipticity condition which ensures that the walk is really
$d$-dimensional. In this paper, the weaker ellipticity Condition~E of~\cite{CP}
will usually not be enough for our purposes. In fact, we believe that most of
our results do not generally hold if one only assumes Condition~E.

For technical reasons we need also the following two conditions:

\medskip
\noindent
{\bf Condition D.} There exists a positive constant~$D_0$ such that the 
expectation of the total number
of the immediate descendants of any particle is at most~$D_0$, i.e.,
$\IP\big[\sum_{v\in\V}|v|\omega_0(v) \leq D_0 \big]=1$.

\medskip
\noindent
{\bf Condition A.} There exist $ x \in \A , v \in \V$ with $\|x\|$
even and $v_x \geq 1$
such that $\IP[\omega_0(v) >0]>0$.

\medskip

We refer to Condition~A as the ``aperiodicity condition'' because, without it,
 the process starting from the origin would live
on even sites at even times, and on odd sites at odd times.
If (with $x,v$ of Condition~A) a site~$y$ is such that
$\omega_y(v)>0$, we say that this site is an \emph{aperiodic}
site, and when $\omega_y(v)>\delta$, we say that this site is 
 $\delta$-aperiodic.
%%   If Condition~A does not hold, then the process (starting from the origin)
%%   lives on even sites (i.e., on the sites for which the sum of coordinates is%%    even)
%%    at even moments and on odd sites at odd moments. It can be obtained 
%%  that the 
We briefly mention that 
suitably adjusted versions of all our results are also true
%% in this case. 
without  Condition~A; the proofs are either the same or even simpler, 
since we
do not have to care about searching for the aperiodic sites, e.g.\ the proofs
of Theorems~\ref{t_shape_reach} and~\ref{t_complete_shape}.

\begin{rem} 
\label{rem:usuel}
{\rm A particular case of the model considered here is the
usual construction of the branching random walk, that was already mentioned
in the beginning of this paper:
for each~$x$,  specify the transition probabilities $p{(x,y)}$,
$y\in \A$, and branching probabilities ${r}_i(x)$, $i=1,2,3,\ldots$.
A particle in~$x$ is first substituted by~$i$ particles with
probability ${r}_i(x)$, then each of the offspring jumps independently
to~$x+y$ with probability $p{(x,y)}$. The pairs $(({
  r}_i(x))_{i\geq 1}, (p{(x,y)})_{y\in \A})$ are chosen
according to some i.i.d.\ field on~$\Z^d$.
In our notations, $\omega_x$ is a mixture of multinomial
distributions on~$\A$:
\[
   \omega_x(\cdot) = \sum_{i\geq 1} {r}_i(x) 
        \mathop{\text{\sf Mult}}(i;p{(x,y)}, y\in\A)(\cdot).
\]
}
\end{rem}

\smallskip
\noindent
{\bf Statement of the results.}
All through the paper, we will assume the above five conditions.
Following~\cite{CP}, we define the notions of transience and recurrence:

\begin{df}
\label{def_rec}
For the particular realization of the random environment~$\om$,
the branching random walk is called recurrent if
\[
\Po^0[\mbox{the origin is visited infinitely often}]=1.
\]
Otherwise, the branching random walk is called
transient.
\end{df}

The following result is a consequence of Propositions~1.2, 1.3, 
and Theorem~1.6 of~\cite{CP}, 
and also Theorem~3.2 of~\cite{M2}; in Section~\ref{s_ind_rw}, we comment about 
the validity of~(i) and~(ii) below in the non-i.i.d.\ case $\rho>1$.

\begin{prop}
\label{rec/tr}
We have either:
\begin{itemize}
\item[(i)] For $\IP$-almost all~$\om$, the branching random walk is recurrent, in
which case $\Po^x[\mbox{the origin is visited infinitely often}]=1$ 
for all~$x\in\Z^d$,
or:
\item[(ii)] For $\IP$-almost all~$\om$, the branching random walk is transient, in
which case $\Po^x[\mbox{the origin is visited infinitely often}]=0$ 
for all~$x\in\Z^d$.
\end{itemize}
The classification criterion: In the case of an i.i.d.\ 
environment, i.e., $\rho=1$,
the \BRW{} in random environment is transient
if and only if there exist $s\in\R^d\setminus\{0\}$ and $\lambda>0$ such that
$\IP$-almost surely we have
\begin{equation}
\label{eq_C_L}
 \sum_{y\in \A} \mu_y^\omega \lambda^{y\cdot s} \leq 1,
\end{equation}
where $\mu_y^\omega=\Eo\eta^0_1(y)$ is the mean number of particles sent to~$y\in \A$ by
one particle at the origin.
\end{prop}
It is interesting to observe that, in the case of an i.i.d.\ 
environment, one has a simple explicit criterion of transience/recurrence
for the \BRW{} in random environment; however, for the many-dimensional
random walk \emph{without branching} in random environment the problem
of finding such a criterion is still far from being solved. In~\cite{CP}
one can find more evidence that branching makes random walks in
random environment ``simpler''; see, for instance, the results 
about the tails of first hitting times.
 
\medskip

Now, we are ready to formulate the main results of this paper. In what follows,
for any $A\subset \R^d$, $\overline A$ is the closure of~$A$, and~$A^o$ is the interior of~$A$.

First, we obtain a shape result about the sites that \emph{can} contain a
particle at time~$n$. 
\begin{theo}
\label{t_shape_reach}
There exists a compact convex set $F\subset \R^d$ such that
for any $\eps>0$, for almost all~$\om$ there exists $n(\om,\eps)$ such that
\begin{equation}
\label{eq_i.1}
\Po[\eta^0_m(x) = 0 \text{ for all } x\in\Z^d\setminus (1+\eps)mF] = 1
\end{equation}
and
\begin{equation}
\label{eq_i.2}
\Po[\eta^0_m(x) = 0] < 1 \qquad \text{ for all } x\in\Z^d\cap (1-\eps)mF
\end{equation}
for all $m\geq n(\om,\eps)$.
\end{theo}

Next, we characterize the (quenched) mean local number of particles.
\begin{theo}
\label{t_local_size}
There exists a convex downwards, continuous function $\beta: F^o \to \R$ 
such that for any closed $G\subset F^o$
\begin{equation}
\label{shape_exp}
\max_{x\in nG\cap\Z^d} |n^{-1}\ln \Eo\eta^0_n(x) - \beta(x/n)| \to 0 \qquad \text{ $\IP$-a.s.},
\end{equation}
 as $n\to\infty$.
\end{theo}
The function $\beta$ is called the local growth exponent for the \BRW.
Denote $B=\{x\in\R^d : \beta(x) \geq 0\}$; by the convexity of~$\beta$, 
the set~$B$ is compact and convex. All through, convex functions will mean
   convex downwards.

\begin{theo}
\label{t_rec_trans}
The \BRW{} in random environment is recurrent if and only if $\beta(0)>0$.
\end{theo}
Note that the border case  $\beta(0)=0$ is transient.
Note also that, from the property~(\ref{sup_beta>0}) below, one obtains that $\beta(0)>0$
if and only if $0\in B^o$.
Thus, given~$B$, one can 
determine whether the \BRW{} is recurrent or transient.

The next result does not only tell us, similarly to Theorem~1.10 of~\cite{CP},
where the particles are
located at time~$n$, but it also gives
an important information about the local size of the population.
\begin{theo}
\label{t_complete_shape}
For any closed $G\subset B^o$,
\begin{equation}
\label{complete_shape_exp}
\max_{x\in nG\cap\Z^d} |n^{-1}\ln \eta^0_n(x) - \beta(x/n)| \to 0
\end{equation}
as $n\to\infty$, a.s.
\end{theo}

\begin{rem} 
\label{rem_complete_shape} 
(Equivalence of models.) {\rm 
The local growth exponent $\beta$ is defined by (\ref{shape_exp}). By 
definition it
only depends on the quenched expectation of the number of particles, and
further, only on the mean number $\mu_y^\omega$
of particles sent in one step to~$y$ by
one particle at the origin. Indeed,
\[
 \Eo\eta^0_n(x)= \sum \prod_{i=1}^n \mu_{x_i-x_{i-1}}^{\omega^{(i)}}
\]
where
$\omega^{(i)}$ is the environment shifted by $x_{i-1}$,
and where the sum ranges over all sequences  $(x_i; 0 \leq i \leq n)$
with $x_0=0, x_n=x,  x_i-x_{i-1}\in \A$. In particular, for any mapping 
$\omega \mapsto \tilde \omega$ from $\M$ to itself such that 
$ \mu_{y}^{\omega}=\mu_{y}^{\tilde \omega}$ for all $y \in \A, \omega \in 
\M$, the two \BRW s in the environments~$\om$ and~$\tilde \om$
are equivalent, in the sense that they have the same local behavior at the
logarithmic scale. Fine details of the branching and displacement
do not matter, under the above five conditions. Moreover, for results
that only concern the expected number of particles (such as Theorem~\ref{t_local_size})
we really do not need the assumption that any particle produces at least one offspring;
this can be substituted by a weaker assumption $\Eo\ZZ^0_1\geq 1$ $\IP$-a.s.}
\end{rem}

Finally, we formulate a result about the total size of the population:
\begin{theo}
\label{t_total_size}
The function~$\beta$ has the property
\begin{equation}
\label{sup_beta>0}
 \sup_{x\in F^o} \beta(x) > 0. 
\end{equation}
For the total size of the population~$\ZZ_n^0$, it holds that
\begin{equation}
\label{eq_total_size}
\frac{\ln \ZZ_n^0}{n} \to \sup_{x\in F^o} \beta(x)
\end{equation}
a.s., as $n\to\infty$.
\end{theo}

Note that, by~(\ref{sup_beta>0}), $B\neq \emptyset$, and $B^o\neq \emptyset$
by convexity. From 
Theorems~\ref{t_local_size} and~\ref{t_complete_shape}
it follows that for $\IP$-almost all~$\om$, for any $\eps>0$
\begin{align*}
\Po[\text{for all $m$ large enough }\eta^0_m(x)=0 \; \forall x\in(1+\eps)mB,~~~&\\
  \eta^0_m(x)\geq 1 \; \forall x\in(1-\eps)mB]&=1.
\end{align*}
As opposed to Theorem~1.10 of~\cite{CP}, here we prove this result 
both in transient and recurrent cases. 

\begin{example} \label{ex1} 
{\rm
{\bf (Constant branching)} More information can be obtained 
in the particular case when
there exists a constant~$\mu>1$ such that
\[
\sum_{y \in \A} \mu^{\omega}_y = \mu \qquad \IP\text{-a.s.} 
\] 
In this case, the expected size of the population is $\Eo \ZZ_n^0 = \mu^n$.
This is the case for instance when $\mu^{\omega}_y$ does not depend on 
$\omega$, and then the \BRW{} is equivalent to a tree-indexed Markov chain
(since in this case we can suppose also that the immediate descendants jump
independently, and the offspring distribution does not depend on the site,
cf.\ Remark~\ref{rem_complete_shape}). 
In the general case, we can define the transitions
\[
p^{\omega_0}(y)= \frac{\mu^\omega_y}{\mu}\,,\quad y \in \A
\]
(recall that $\mu^\omega_y$ only depends on $\omega_0$), 
and consider the random walk in random environment~$\chi_n$ with transition
probabilities
$\Po[\chi_{n+1}=x+y \mid \chi_n=x]=p^{\omega_x}(y)$ for $y \in \A$.  
We see here that
$ \Eo\eta_n^0(x)=\mu^n \Po(\chi_n=x)$, and therefore, by Theorem~\ref{t_local_size},
\[
\beta(a)=\ln \mu + \lim_{n \to \8 } n^{-1} \ln \Po\big[ \chi_n = [na]\big] 
\]
with~$[na]$ the integer part of~$na$ (coordinatewise). The limit can be 
expressed in terms of the quenched large deviation rate function~$I^q$,
which have been studied in the nestling case in \cite{Zerner} and in complete 
generality in~\cite{Varadhan}: it holds that
\[
\beta(a) = \ln \mu - I^q(a).
\]
These references are for i.i.d.\ environment and estimate the probability 
of sets in the scale~$n$ instead of the probability of points, but 
one can see that they apply to our discussion here. This example shows that 
the convex function $\beta$ is not necessarily strictly  convex. Indeed,
it is known that the rate function has a flat horizontal part if the
random walk in random environment is nestling with a non-zero speed $v$, in
which case   $I^q(a)=0$ for $a$ in the whole 
interval with endpoints $0$ and $v$. Finally, it is straightforward to see
from the nearest neighbor jumps case, that the shapes $F$ and $B$ may have
``facets'', i.e., flat parts on their boundaries (cf.\ e.g.\ Example~7 of~\cite{CP}).  
}
\end{example}

\section{Proofs}
\label{s_proofs}
The rest of this paper is organised as follows. First, in Section~\ref{s_ind_rw} we recall some
concepts of~\cite{CP}, such as recurrent seeds and induced random walks. In Section~\ref{s_proof_shape_reach}
we study the set of sites which can be reached up to time~$n$ (this amounts, basically, to applying
the Subadditive Ergodic Theorem). In Section~\ref{s_proof_t_local_size} we prove the results related
to the \emph{expected} local population size, and in Section~\ref{s_proof_t_complete_shape}
we study the local population size itself. Finally, in Section~\ref{s_proof_exp_growth} we
prove the equation~(\ref{sup_beta>0}), i.e., that the expectation of the total number of particles
grows exponentially 
(somewhat surprisingly, this is one of the most difficult results of this paper),
and then we prove Theorem~\ref{t_total_size}.

\subsection{Induced random walks}
\label{s_ind_rw}
To begin, we introduce some more basic notations.
% :
% for $x=(x^{(1)},\ldots,x^{(d)})\in\Z^d$ write
% \[
%  \|x\|_{\infty} = \max_{i=1,\ldots,d} |x^{(i)}|.
% \]
We denote by~$\Q$ the set of rational numbers, and define $\N:=\{1,2,3,\ldots\}$.
Let~$L_0$ to be the maximal jump length, i.e.,
\[
 L_0 = \max_{x\in \A} \|x\|,
\]
and let~$\K_n$ be the $d$-dimensional discrete ball with respect to 
the 
%%${\mathcal L}_1$-norm:
$\ell_1$-norm:
\begin{equation}
\label{def_Kn}
 \K_n = \{x\in\Z^d: \|x\|\leq n\}.
\end{equation}

As in Section~2.1 of~\cite{CP}, we define now the notion of \emph{induced} random walk
in random environment associated with
the branching random walk in  random environment.
Defining
\[
\tV= \big\{ (v,\kappa): v \in \V, \kappa {\rm \ probability \ measure\ on \ }
\{y: v(y)\geq 1\} \big\},
\]
we consider some probability measure $\tilde\IP$ on $\tV^{\Z^d}$ with marginal~$\IP$
on $\V^{\Z^d}$. A stationary random field
$\tilde \om = ((\omega_x, \kappa_x), x \in \Z^d)$  with the law $\tilde\IP$
defines our branching random walk as above, coupled with
a random walk in random environment with transition
probability
\begin{equation}
\label{prob_ind_rw}
p_x(y)= \sum_{v \in \V} \omega_x(v) \kappa_x(y)
\end{equation}
from $x$ to $x+y$. In words, we pick randomly one of the children in
the branching random walk. To keep things simple, we will drop the 
tilde from the notations $\tilde\IP$.
In this paper, we need only the so-called uniform induced random walk,
for which the measure~$\kappa$ is defined as follows: 
$\kappa$ is uniform on the locations $\{x \in \A: v_x\geq 1\}$.

% \begin{itemize}
% \item[(i)] uniform: $\kappa$ is uniform on the locations
% $\{x \in \A: v_x\geq 1\}$;
% \item[(ii)] particle-uniform: $\kappa(y)$ is proportional to the number of
% particles sent by~$v$ to~$y$.
% \end{itemize}
% 
% {\bf Do we need to introduce the 2 of them ?? }

An important idea that will be repeatedly employed in this paper is 
to use the uniform ellipticity of the walk in order to reveal some
 independence in the environment:
because of Condition~UE, the uniform induced random walk is uniformly 
elliptic
as well, 
% \footnote{ Under Condition~D, the particle-uniform induced random walk
% is uniformly elliptic too.}, 
 and so sometimes it makes its steps ``without looking at the 
environment''.
A similar construction can be found in ~\cite{CZ,Z}.
%% e.g.~\cite{CZ,Z}).
Specifically, let us consider the uniform induced random walk $\xi^z$
($z$ stands for the starting location of this random
walk). According to~(\ref{prob_ind_rw}), the transition probabilities for $\xi^z$ are:
\[
 \sigma(x,x+y) = \Po[\xi^z_{n+1}=x+y \mid \xi^z_n=x] 
       = \sum_{v:v_y\geq 1}\frac{\omega_x(v)}{|\{u:v_u\geq 1\}|},
\]
which means that, in the case when a particle has more than one offspring
in the branching random walk (i.e., it produces a configuration~$v$ with $|v|>1$), 
the next (relative) location for the uniform induced random walk
is chosen uniformly among the locations $\{x \in \A: v_x\geq 1\}$.
By Condition~UE, this induced random walk is uniformly elliptic
in the sense that 
\[
 \Po[\xi^z_{n+1}=x+e \mid \xi^z_n=x]\geq {\hat\eps}_0 
           \qquad \text{for all }e\in\{\pm e_i, i=1,\ldots,d\}
\]
where ${\hat\eps}_0=\eps_0/|\A|$.
%%    Let $Z_1,Z_2,Z_3,\ldots$ be a sequence of i.i.d.\ random variables,
%%    $\Po[Z_i=1]=1-\Po[Z_i=0]=2d{\hat\eps}_0$. Let us enumerate the
%%    elements of the set $\{\pm e_i, i=1,\ldots,d\}$ in some order, so that
%%    $\{\pm e_i, i=1,\ldots,d\} = \{{\hat e}_i, i=1,\ldots,2d\}$.
%%    The construction of the induced random walk can now be performed in
%%    the following way:
%%    \begin{itemize}
%%   \item if $Z_i=1$, then $\xi^z_i = \xi^z_{i-1}+{\hat e}_j$ with probability
%%    $(2d)^{-1}$, $j=1,\ldots,2d$;
%%   \item if $Z_i=0$, then $\xi^z_i = \xi^z_{i-1}+{\hat e}_j$ with probability
%%   $\frac{\sigma(\xi^z_{i-1},\xi^z_{i-1}+{\hat e}_j)-
%%    {\hat\eps}_0}{1-2d{\hat\eps}_0}$,
%%    and $\xi^z_i = \xi^z_{i-1}+y$ with probability
%%     $\frac{\sigma(\xi^z_{i-1},\xi^z_{i-1}+y)}{1-2d{\hat\eps}_0}$ for
%%    $y\notin \{{\hat e}_i, i=1,\ldots,2d\}$.
%%    \end{itemize}
Let $\hat Z_1,\hat Z_2,\hat Z_3,\ldots$ be a sequence of i.i.d.\ 
random variables with values in $\{0,1,\ldots,2d\}$, such that
$\Po[\hat Z_i=j]={\hat\eps}_0$ for $j \neq 0$ and of course $\Po[\hat Z_i=0]=
1-2d{\hat\eps}_0$. We still keep the symbol $\Po$ to denote the probability 
on the enlarged 
probability space where both the \BRW~and the  sequence $(\hat Z_i)_i$ are 
defined. Set $Z_i=\1{\hat Z_i \neq 0}  $.
Let us enumerate the
elements of the set $\{\pm e_i, i=1,\ldots,d\}$ in some order, so that
$\{\pm e_i, i=1,\ldots,d\} = \{{\hat e}_i, i=1,\ldots,2d\}$.
We can now construct the induced random walk as 
follows:
\begin{itemize}
 \item if $\hat Z_i=j$ for some $j \neq 0$, then 
$\xi^z_i = \xi^z_{i-1}+{\hat e}_j$;
 \item if $\hat Z_i=0$, then $\xi^z_i = \xi^z_{i-1}+{\hat e}_j$ 
with probability
$\frac{\sigma(\xi^z_{i-1},\xi^z_{i-1}+{\hat e}_j)-{\hat\eps}_0}
{1-2d{\hat\eps}_0}$,
and $\xi^z_i = \xi^z_{i-1}+y$ with probability
 $\frac{\sigma(\xi^z_{i-1},\xi^z_{i-1}+y)}{1-2d{\hat\eps}_0}$ for
$y\notin \{{\hat e}_i, i=1,\ldots,2d\}$.
\end{itemize}
In words, this means that when the value of $Z$-variable is~$1$, the random
walk moves without looking at the random environment.

In the next definition we recall the notion of $(U,\HH)$-{\em seed} (suitably
adapted for the case of finitely dependent environment), that was introduced in~\cite{CP}.

\begin{df}
\label{def_seed}
Fix a finite set $U \subset \Z^d$ containing 0,
and $\HH_x \subset \M$ with $\IP[\omega_x\in \HH_x \text{ for all }x \in U]>0$. 
With $\HH=(\HH_x, x \in U)$, the pair $(U,\HH)$ is called a seed.
We say that $\om$ has a $(U,\HH)$-seed at $z \in \Z^d$
(or that a $(U,\HH)$-seed occurs in~$z$) if
\[
\omega_{z+x} \in \HH_x  \mbox{ ~for all~ } x \in U,
\]
and that~$\om$ has a $(U,\HH)$-seed in the case $z=0$.
We call~$z$ the center of the seed.
\end{df}

As in~\cite{CP}, it is easy to see that with probability~$1$
the branching random walk visits infinitely many distinct $(U,\HH)$-seeds
(this can be done by showing that the uniform induced random walk does so). 
We now give the argument. With $r=\rho + {\rm diameter}(U)$,
at any time a.s.\ there exist subsequent times $n, t=n+r$
such that at time~$t$ the uniform induced random walk is situated in 
 a location~$x$ which is at distance~$r$ away from its range up to time~$n$,
without looking at the environment 
(that is, $Z_{n+1}=\ldots=Z_t=1$). Then, by
Condition~I, the environment can be constructed inside the translate $x+U$ 
independently from all what done before, and so the probability to 
generate a $(U,\HH)$-seed at site~$x$ is a positive constant. 
By the Borel-Cantelli lemma, with probability~1 an infinite number of 
$(U,\HH)$-seeds will be visited.

Then, still following~\cite{CP}, we define the branching random walk 
\emph{restricted} 
to set~$M\subset\Z^d$ simply by
discarding all particles that step outside~$M$, and write~$\PoA{M}, \EoA{M}$
for corresponding probability and expectation. Next, we consider a
shortened version of Definition~2.5 from~\cite{CP}:
\begin{df}
\label{rec_seed}
Let~$U$ be a finite subset of~$\Z^d$ with $0 \in U$.
Let~$\p$ be a probability distribution on $\Z_+$
 with mean larger than 1, i.e.,
$\p=(p_0, p_1,p_2,\ldots)$ with $p_i\geq 0$,
$\sum p_i=1$, $\sum ip_i>1$.
An $(U,\HH)$-seed is called $\p$-recurrent if for any~$\om$
such that $\omega_x\in \HH_x, x\in U$, we have
\[
\PoA{U}^0[ 0 \text{ will be visited by at least } i \text{ ``free" particles}]
   \geq \sum_{j=i}^\infty p_j
\]
for all $i \geq 1$. By ``{\em free}'' particles
we mean that none is the descendant of
another one. 
\end{df}
It is important to note that, by definition of the restricted branching random walk,
the above probability
only depends on the environment inside the $\rho$-neighborhood of~$U$.

Then, it is straightforward to see that all the discussion of Section~2.2
of~\cite{CP} readily extends to the case of finitely dependent 
environment as well. In particular, the recurrence is equivalent to the 
existence of recurrent seeds; this fact will be used several times in this paper.

\subsection{Proof of Theorem~\ref{t_shape_reach}}
\label{s_proof_shape_reach}
For arbitrary $x,y\in \Z^d$, $\delta\in [0,\eps_0]$ ($\eps_0$ is
from Condition~UE), define
\begin{align*}
 T_\omega^\delta(x,y) &= \min\{n: \text{ there exist }z_0,z_1,\ldots,z_n\in\Z^d
   \text{ with }z_0=x, z_n=y  \\
    & \qquad \qquad \text{ such that }\omega_{z_i}(v:v_{z_{i+1}-z_i}\geq 1)>\delta, \;\;
          i=0,\ldots,n-1\},
\end{align*}
so that $T_\omega^\delta(x,y)$ is the minimal number of steps necessary for a particle
in~$x$ to send an offspring to~$y$, with the condition also that this should happen
with big enough probability on each step. By Condition~UE, it is immediate that
\begin{equation}
\label{T<||_1}
 T_\omega^\delta(x,y) \leq \|x-y\|.
\end{equation}
Clearly, this family of random variables has the subadditive property:
for any $x,y,z\in\Z^d$, $\delta\in[0,\eps_0]$, and \emph{any}~$\om$
\begin{equation}
\label{T_subadd}
 T_\omega^\delta(x,y) \leq T_\omega^\delta(x,z) + T_\omega^\delta(z,y).
\end{equation}

Consider any $a\in \Q^d$ and define
\begin{equation}
\label{eq_def_mu}
 \mu^\delta (a) = \lim_{n\to\infty} \frac{T_\omega^\delta(0,k_0 an)}{k_0 n},
\end{equation}
where $k_0$ is the smallest positive integer such that $k_0a\in\Z^d$.
With~(\ref{T<||_1}) and~(\ref{T_subadd}), the Subadditive Ergodic Theorem
(see e.g.\ Theorem~2.6 of Chapter~VI of~\cite{L_book}) 
shows that the (nonrandom) limit
in~(\ref{eq_def_mu}) exists a.s.\ and in ${\mathcal L}_1$;
by~(\ref{T_subadd}), this limit verifies
$\mu^\delta(a)+\mu^\delta(b)\geq\mu^\delta(a+b)$,
$\mu^\delta(ra)=r\mu^\delta(a)$, for any
$a,b\in\Q^d, r\in\Q^+$. Moreover, since the jumps are bounded, we have
$ T_\omega^\delta(x,y) \geq L_0^{-1}  \|x-y\|$ and finally
$L_0^{-1}  \|a\|\leq \mu^\delta(a) \leq \|a\|$.
Then, by continuity
one can define~$\mu^\delta(a)$ for any~$a\in\R^d$ in such a way 
that~$\mu^\delta$ is a norm on~$\R^d$. 

Let 
\[
 F_\delta = \{a\in\R^d : \mu^\delta(a)\leq 1\}.
\]
Clearly,
 for any $\delta\in [0,\eps_0]$, the set $F_\delta$ is compact and
convex, and $0\in F_\delta^o$. 
By definition, for any $x,y$ and~$\om$ it holds that
\begin{equation}
\label{monot_delta}
 T_\omega^{\delta_1}(x,y) \leq T_\omega^{\delta_2}(x,y)
\end{equation}
when $\delta_1\leq \delta_2$.
So, we have that $F_{\delta_1}\subset F_{\delta_2}$
for $\delta_1\geq \delta_2$. For the rest of this paper, denote $F:=F_0$
(and this is the compact convex set we are looking for in Theorem~\ref{t_shape_reach}).
The following lemma shows that the family $F_\delta$ is continuous in~$\delta$:
\begin{lmm}
\label{l_cont_delta}
%%%%%%%%%%%%  NON !!! 
%% We have $\lim_{\delta\to 0} \dist (F_\delta, F)=0$.
%% cela ne correspond pas a la def de \dist donnee au debut
For any $\eps>0$ there exists $\delta>0$ such that
$(1-\eps)F \subset F_\delta \subset F$.
\end{lmm}

\noindent
{\it Proof: }
%%Clearly, 
By covering $F$ with finitely many small disks, we see 
it is enough to prove that, for any $a\in\Q^d$,
\begin{equation}
\label{eq_cont_delta}
 \lim_{\delta\to 0} \Big[\lim_{n\to\infty}\frac{T_\omega^\delta(0,k_0 an)}{k_0 n}
        - \lim_{n\to\infty}\frac{T_\omega^0(0,k_0 an)}{k_0 n}\Big] = 0,
\end{equation}
where $k_0$ is the smallest positive integer such that $k_0a\in\Z^d$.
Observe that, by the Subadditive Ergodic Theorem, the left-hand side
of~(\ref{eq_cont_delta}) is equal to
\[
 \lim_{\delta\to 0} \Big[\inf_{n\geq 1}\frac{\IE T_\omega^\delta(0,k_0 an)}{k_0 n}
        - \inf_{n\geq 1}\frac{\IE T_\omega^0(0,k_0 an)}{k_0 n}\Big].
\]
 Note that, by~(\ref{T<||_1}), 
$T_\omega^\delta(x,y)$ depends only on a finite piece of the environment,
so
\begin{equation}
\label{shodimost'_T}
 \IE T_\omega^\delta(x,y) \to \IE T_\omega^0(x,y)
\end{equation}
as $\delta\to 0$. Now, fix an arbitrary $\eps>0$ and 
choose $n_1$ in such a way that 
\[
 \frac{\IE T_\omega^0(0,k_0 an_1)}{k_0 n_1} - 
      \inf_{n\geq 1}\frac{\IE T_\omega^0(0,k_0 an)}{k_0 n} < \eps.
\]
By~(\ref{monot_delta}) and~(\ref{shodimost'_T}), there exists~$\delta_1>0$
such that 
\[
 \IE T_\omega^\delta(0,k_0 an_1) - \IE T_\omega^0(0,k_0 an_1) < k_0n_1\eps
\]
for all $\delta\leq \delta_1$. So, we obtain that
\[
 \limsup_{\delta\to 0} 
 \Big[\inf_{n\geq 1}\frac{\IE T_\omega^\delta(0,k_0 an)}{k_0 n}
        - \inf_{n\geq 1}\frac{\IE T_\omega^0(0,k_0 an)}{k_0 n}\Big] \leq 2\eps,
\]
which implies~(\ref{eq_cont_delta}).
\qed

Now, we are ready to prove the first part of Theorem~\ref{t_shape_reach}.
Denote
\[
 W_\omega^\delta(n) = \{x\in\Z^d : T_\omega^\delta(0,x)\leq n\},
\]
and let ${\hat W}_\omega^\delta(n) = W_\omega^\delta(n) + (-1/2,1/2]^d$.

By a standard argument (see e.g.~\cite{AMP,BG,DG}) 
one can show that, for any $\eps>0$,
\begin{equation}
\label{shape_t_W} 
 (1-\eps)F_\delta \subset \frac{{\hat W}_\omega^\delta(n)}{n}
 \subset (1+\eps)F_\delta
\end{equation}
for all~$n$ large enough.
In particular, since $n^{-1}{\hat W}_\omega^\delta(n)\subset (1+\eps)F$
for all~$n$ large enough, the first claim of Theorem~\ref{t_shape_reach}
follows.

In order to prove the second claim, let us define by
\begin{align}
 R_\omega^{x,\delta}(n) &= \{y\in\Z^d : 
          \text{ there exist }z_0,z_1,\ldots,z_n\in\Z^d
   \text{ with }z_0=x, z_n=y  \nonumber\\
    & \qquad \text{ such that }\omega_{z_i}(v:v_{z_{i+1}-z_i}\geq 1)>\delta, \;\;
          i=0,\ldots,n-1\}, 
\label{def_Rn}
\end{align}
the set of sites that can be reached in \emph{exactly}~$n$ steps
(with our usual restriction on the probabilities of the steps).
Clearly, if $y\in R_\omega^{x,\delta}(n)$, then $\Eo\eta^x_n(y)\geq\delta^n$.
Denoting also ${\hat R}_\omega^{x,\delta}(n)=R_\omega^{x,\delta}(n)+(-1/2,1/2]^d$,
we intend to prove that, for any~$\eps>0$ and almost all~$\om$
\begin{equation}
\label{incl_Rn}
 (1-\eps)F_\delta \subset \frac{{\hat R}_\omega^{0,\delta}(n)}{n}
 \subset (1+\eps)F_\delta
\end{equation}
for all~$n$ large enough.

Let (recall the definition of $\K$ from~(\ref{def_Kn}))
\begin{align*}
 M_n^\delta &= 
\{\om : \text{ for any }x\in \K_{L_0n} \text{ there exists a $\delta$-aperiodic
site }y \\
 & \qquad \qquad\text{ such that } \|x-y\|\leq n^{1/2}\}.
\end{align*}
Since the random environment is finitely dependent, one obtains that there 
are some
positive constants $\delta_0,C_1,C_2$ such that for all~$\delta\leq \delta_0$
\begin{equation}
\label{ocenka_Mn}
 \IP[M_n^\delta] \geq 1-C_1n^d\exp(-C_2 n^{-d/2}).
\end{equation}
By Borel-Cantelli lemma,
\begin{equation}
\label{BC_Mn}
 \IP[\text{there exists $n(\om)$ such that }M_n^\delta \text{ occurs
       for all }n\geq n(\om)]=1.
\end{equation}

Fix any $\eps>0$ and consider a site $x\in (1-\eps)nF_\delta$.
As we know from~(\ref{shape_t_W}), if~$n$ is large enough, then 
$T_\omega^\delta(0,x)\leq (1-\frac{\eps}{2})n$ and the event~$M_n^\delta$ occurs.
Now, consider two cases:
\begin{enumerate}
 \item $n-T_\omega^\delta(0,x)$ is even. Then it is trivial to
obtain that $x\in R_\omega^{0,\delta}(n)$ (one can complete the path
of length $T_\omega^\delta(0,x)$ which ends in~$x$ by $x+e_1,x,x+e_1,x,\ldots$).
 \item $n-T_\omega^\delta(0,x)$ is odd.
Suppose also that~$n$ is so large that 
$\max\{L_0,n^{1/2}\}<\frac{\eps n}{6}$. 
Since~$M_n^\delta$ occurs, there exists an aperiodic site~$y$ such that
$\|x-y\|<\frac{\eps n}{6}$. Then, to complete the path
of length $T_\omega^\delta(0,x)$ which ends in~$x$, essentially
one goes from~$x$ to~$y$ in exactly $\|x-y\|$ steps,
then jumps from~$y$ to some~$y_1$ with $\|y-y_1\|$ even, 
then goes back to~$x$ (and then, if necessary, one puts
$x+e_1,x,x+e_1,x,\ldots$ to the end of the path).
\end{enumerate}

In both cases we obtain that $x\in R_\omega^{0,\delta}(n)$,
and this concludes the proof of~(\ref{incl_Rn}) 
and thus of Theorem~\ref{t_shape_reach}.
\qed

\subsection{Proof of Theorems~\ref{t_local_size} and~\ref{t_rec_trans}}
\label{s_proof_t_local_size}
We begin by showing that the function~$\beta$ can be defined in the following way:

\begin{lmm}
\label{l_def_beta}
For any $a\in F^o\cap \Q^d$, $a\neq 0$, the following
quantity $\beta(a)$ is well-defined and is a.s.\ constant:
\begin{equation}
\label{eq_def_beta}
 \beta(a) = \lim_{n\to\infty} 
\frac{\ln\Eo\eta^0_{k_0n}(k_0na)}{k_0n} \qquad\text{a.s.},
\end{equation}
where $k_0$ is the smallest positive even integer number such that $k_0a\in 2\Z^d$.
\end{lmm}

\noindent
{\it Proof: } The expected number of particles $\Eo\eta$
has a supermultiplicative property: for any $x,y,z\in\Z^d$, $n_1,n_2\geq 0$
\begin{equation}
\label{supermult}
 \Eo\eta^x_{n_1}(y)\Eo\eta^y_{n_2}(z) \leq \Eo\eta^x_{n_1+n_2}(z),
\end{equation}
so the family of random variables 
\[
 S_{m,n}=k_0^{-1}\ln \Eo\eta^{k_0ma}_{k_0(n-m)}(k_0na)
\]
 is superadditive.
Note, however, that the random variables of the
latter family may assume the value $-\infty$.

Suppose first that $\|a\|\leq 1$. Then, from Condition~UE
we obtain that
\[
 \Eo\eta^0_{k_0}(k_0a) \geq \eps_0^{k_0},
\]
and so, taking Condition~D into account,
 the existence of the limit in~(\ref{eq_def_beta}) immediately
follows from the Subadditive Ergodic Theorem.

However, if one wants to apply the Subadditive Ergodic Theorem
to the family $S_{m,n}$ in the case $\|a\|>1$, there is the following difficulty: 
there may be some $a\in F\cap \Q^d$, such that with positive 
probability it happens that $T_\omega^0(0,k_0na)>k_0n$,
which means that $\IE(-\ln\Eo\eta^0_{k_0n}(k_0na))^+=\infty$
(even $\IP[(-\ln\Eo\eta^0_{k_0n}(k_0na))^+=\infty]>0$).
So, for the case $\|a\|>1$ we need a different approach. 

For the rest of the proof of Lemma~\ref{l_def_beta} we
suppose that $\|a\|>1$ and let~$\delta$ be such that $a\in F_\delta^o$.
Define
\begin{equation}
\label{def_beta_a}
\beta(a) = \limsup_{n\to\infty} \frac{\ln\Eo\eta^0_{k_0n}(k_0an)}{k_0n}
\end{equation}
(in principle, $\beta(a)$ could depend also on~$\om$, but in the next few lines we
will show that it is a.s.\ constant).
By~(\ref{incl_Rn}), $\IP$-a.s.\ there exists~$n(\om)$ such that (recall~(\ref{def_Rn}))
$k_0an\in R_\omega^{0,\delta}(k_0n)$ for all $n\geq n(\om)$.
Using~(\ref{supermult}), we obtain that for all $n\geq n(\om)$
\[
 \ln\Eo\eta^0_{k_0(n+m)}(k_0a(n+m)) \geq \ln\Eo\eta^0_{k_0n}(k_0an) 
                  + \ln\Eo\eta^{k_0an}_{k_0m}(k_0a(n+m))
\]
for all~$m$ such that $k_0a(n+m)\in R_\omega^{k_0an,\delta}(k_0m)$, which means that
\[
 \limsup_{m\to\infty} \frac{\ln\Eo\eta^0_{k_0m}(k_0am)}{k_0m} \geq
          \limsup_{m\to\infty} \frac{\ln\Eo\eta^{k_0an}_{k_0m}(k_0a(n+m))}{k_0m}
\]
for all $n\geq n(\om)$. Since the sequence
\[
 \Big(\limsup_{m\to\infty} \frac{\ln\Eo\eta^{k_0an}_{k_0m}(k_0a(n+m))}{k_0m}\Big)_{n=0,1,2,\ldots}
\]
is stationary ergodic, this shows that
the upper limit in~(\ref{def_beta_a}) is a.s.\ constant.

Now, our goal is to prove that
\begin{equation}
\label{liminf_beta_a}
\liminf_{n\to\infty} \frac{\ln\Eo\eta^0_{k_0n}(k_0an)}{k_0n} \geq \beta(a).
\end{equation}
Choose $r\in \Q\cap(1,+\infty)$
in such a way that $ra\in F_\delta^o$;
let~$k_1$ be the smallest positive integer
such that $k_1r\in\Z$. Fix a small $\alpha>0$. 
By~(\ref{incl_Rn}), for all~$n$ large enough we
have $k_0k_1ra\lf\alpha n\rf \in R_\omega^{0,\delta}(k_0k_1\lf\alpha n\rf )$.
Recall that, if $y\in R_\omega^{x,\delta}(n)$, then
$\Eo\eta^x_y(n)\geq \delta^n$, so
\begin{equation}
\label{pervyj_shag}
\Eo\eta^0_{k_0k_1\lf\alpha n\rf }(k_0k_1ra\lf\alpha n\rf ) \geq \delta^{k_0k_1\lf\alpha n\rf }.
\end{equation}

To proceed, we use the approach of~\cite{S}. Fix any $\eps>0$ and define
the events
\begin{align*}
 H^m(N) &= \{S_{m,m+k}<k(\beta(a)-\eps) \text{ for all }k=1,\ldots,N\},\\
 G^m(N) &= (H^m(N))^c.
\end{align*}
By definition, we have that
\begin{equation}
\label{prob_H_to_0}
\IP[H^m(N)] \to 0 \qquad \text{as }N\to\infty.
\end{equation}

Now, we divide the integer interval $[k_1r\lf\alpha n\rf ,n)$ into
some subintervals and some singletons using the following
algorithm. Begin with $k=k_1r\lf\alpha n\rf $; 
inductively, let~$k$ be the smallest
integer not yet assigned. If the event $G^k(N)$ occurs, then
there exists $\ell\leq N$ such that $S_{k,k+\ell}\geq \ell (\beta(a)-\eps)$.
In this case we add the interval $[k,k+\ell)$ to our collection
(and then pass to $k'=k+\ell$). On the other hand, if the event 
$H^k(N)$ occurs, then we declare~$k$ to be a singleton 
(and then pass to $k'=k+1$). As a result of this procedure,
we obtain~$u$ intervals $[\tau_i,\tau_i+\ell_i)$, $i=1,\ldots,u$, and~$w$
singletons $\sigma_1,\ldots,\sigma_w$. For each of the above intervals, we
have $S_{\tau_i,\tau_i+\ell_i}\geq \ell_i (\beta(a)-\eps)$, so that
\begin{equation}
\label{good_interval}
 \Eo\eta^{k_0a\tau_i}_{k_0\ell_i}(k_0a(\tau_i+\ell_i)) 
            \geq \exp\big(k_0\ell_i(\beta(a)-\eps)\big),
\end{equation}
$i=1,\ldots,u$.

Then, Condition~UE implies that
\begin{equation}
\label{reach_small_steps}
 \Eo\eta^{k_0am}_{k_0\|a\|}(k_0a(m+1)) \geq \eps_0^{k_0\|a\|}
\end{equation}
for any~$m$. So, abbreviating $\phi_1(n) = \sum_{i=1}^u\ell_i$,
$\phi_2(n) = w = n-k_1ra\lf\alpha n\rf  - \phi_1(n)$, $t_n = k_0\phi_1(n) + k_0\|a\|\phi_2(n)$,
we obtain from~(\ref{supermult}), (\ref{pervyj_shag}), 
(\ref{good_interval}), and~(\ref{reach_small_steps}) that
\begin{align}
 \Eo\eta^{0}_{t_n+k_0k_1\lf\alpha n\rf }(k_0 an) &\geq  
             \Eo\eta^0_{k_0k_1\lf\alpha n\rf }(k_0k_1ra\lf\alpha n\rf )
 \prod_{i=1}^u \Eo\eta^{k_0a\tau_i}_{k_0\ell_i} (k_0a(\tau_i+\ell_i))   \nonumber\\  
 & \qquad {}\times \prod_{j=1}^w \Eo\eta^{k_0a\sigma_i}_{k_0\|a\|}(k_0a(\sigma_i+1)) \nonumber\\
 &\geq  \delta^{k_0k_1\lf\alpha n\rf } 
          \exp\big(k_0\phi_1(n)(\beta(a)-\eps)\big) \eps_0^{k_0\|a\|\phi_2(n)}.
\label{interm_1}
\end{align}

By construction of the intervals, we have
\[
 \sum_{i=1}^u\ell_i \geq n-k_1r\lf\alpha n\rf -N 
               - \sum_{j=k_1r\lf\alpha n\rf }^n \I_{H^j(N)},
\]
so, by Birkhoff's theorem
\begin{equation}
\label{oc_phi1}
 \liminf_{n\to\infty}\frac{\phi_1(n)}{n} \geq (1-\IP[H^0(N)])(1-k_1r\alpha);
\end{equation}
then,
\begin{equation}
\label{oc_phi2}
 \limsup_{n\to\infty}\frac{\phi_2(n)}{n} \leq \IP[H^0(N)](1-k_1r\alpha).
\end{equation}
Take $\alpha$ small so that $k_1r\alpha < 1$.
Note that $t_n=k_0\big(n-k_1r\lf\alpha n\rf +(\|a\|-1)\phi_2(n)\big)$,
so, by~(\ref{oc_phi2}) and~(\ref{prob_H_to_0}),
\begin{equation}
\label{limsup_tn}
 \limsup_{n\to\infty} \frac{t_n}{k_0n} \leq
 (1-k_1r\alpha)\big(1+(\|a\|-1)\IP[H^0(N)]\big)
%%  1-k_1r\alpha\big(1-(\|a\|-1)\IP[H^0(N)]\big)
                 < 1-k_1\alpha
\end{equation}
if~$N$ is so large that 
%%$r\big(1-(\|a\|-1)\IP[H^0(N)]\big)^{-1}>1$.
$\IP[H^0(N)] \leq k_1(r-1)\alpha (1-k_1r\alpha)^{-1} (\|a\|-1)^{-1}$.
Thus, using~(\ref{interm_1}), we obtain
\begin{align*}
 \Eo\eta^{0}_{k_0n}(k_0 an) &\geq  \Eo\eta^{0}_{t_n+k_0k_1\lf\alpha n\rf }(k_0 an)
                  \Eo\eta^{k_0 an}_{k_0n-t_n-k_0k_1\lf\alpha n\rf }(k_0 an)\\
  &\geq  \delta^{k_0k_1\lf\alpha n\rf } 
          \exp\big(k_0\phi_1(n)(\beta(a)-\eps)\big) 
               \eps_0^{k_0(\|a\|\phi_2(n) + n-\frac{t_n}{k_0}-k_1\lf\alpha n\rf )}.
\end{align*}
The inequality~(\ref{liminf_beta_a}) now follows from~(\ref{prob_H_to_0}),
(\ref{oc_phi1}) (note also that, trivially, $\limsup_{n\to\infty}n^{-1}\phi_1(n)\leq 1$),
and~(\ref{oc_phi2}). Letting $\alpha \searrow 0$, this concludes the proof of 
Lemma~\ref{l_def_beta}.
\qed

\medskip

\begin{lmm}
\label{l_convex_beta}
The function $\beta(a)$ is convex downwards on $F^o\cap \Q^d$ (and so it can
be defined for all $a\in F^o$ by continuity, preserving the convexity).
\end{lmm}

\noindent
{\it Proof: }
Consider $a,b\in F^o$
such that $\|a-b\|\leq 1$, and note that
there exists $\delta>0$ such that $a,b\in F^o_\delta$. 
Now, we have to prove that for any $s\in (0,1)\cap \Q$,
\begin{equation}
\label{eq_conv}
\beta(sa+(1-s)b)\leq s\beta(a)+(1-s)\beta(b).
\end{equation}
Let $k_0=\min\{k\in 2\N : k a\in 2\Z^d\}$,
$k_1=\min\{k\in 2\N : k b\in 2\Z^d\}$, and~$\ell=\min\{k\in\N: ks\in\N\}$.
Use the abbreviation
\[
 A^m_n(x) = \1{k_0k_1\ell xn \in R_\omega^{k_0k_1\ell xm,\delta}(k_0k_1\ell(n-m))}.
\]
Note that, by~(\ref{incl_Rn}), Lemma~\ref{l_def_beta},
 and the bounded convergence theorem, for any~$x\in F^o_\delta$
it holds that (with $k'=\min\{k\in 2\N : k x\in 2\Z^d\}$; we use the convention $0\times\infty=0$)
\begin{equation}
\label{bounded_conv}
 \frac{\IE\big(\1{k'xn\in R_\omega^{0,\delta}}\ln\Eo\eta^0_{k'n}(k'xn)\big)}{k'n} \to \beta(x) 
     \qquad\text{as $n\to\infty$}.
\end{equation}

 From the supermultiplicative property~(\ref{supermult}) we obtain
\begin{align}
\lefteqn{\frac{\ln\Eo\eta^0_{k_0k_1\ell n}(k_0k_1\ell (sa+(1-s)b) n)}{k_0k_1\ell n}
 \I_{A^0_{sn}(sa)}\I_{A^{sn}_n(sa+(1-s)b)}}
\nonumber\\
 &\geq
 s\frac{\ln\Eo\eta^0_{k_0k_1\ell s n}(k_0k_1\ell sa n)}{k_0k_1\ell sn}
 \I_{A^0_{sn}(sa)}\I_{A^{sn}_n(sa+(1-s)b)} \label{eq_conv_proof}\\
  & \quad {}+ (1-s)\frac{\ln\Eo\eta^{k_0k_1\ell sa n}_{k_0k_1\ell (1-s) n}
         (k_0k_1\ell (sa+(1-s)b) n)}{k_0k_1\ell (1-s)n}\I_{A^0_{sn}(sa)}\I_{A^{sn}_n(sa+(1-s)b)}.
\nonumber
\end{align}
Since 
\begin{align*}
\lefteqn{\I_{A^{sn}_n(sa+(1-s)b)}\Eo\eta^{k_0k_1\ell sa n}_{k_0k_1\ell (1-s) n}
(k_0k_1\ell (sa+(1-s)b) n)}\qquad\\
& \eqlaw 
        \I_{A^0_{(1-s)n}((1-s)b)}\Eo\eta^0_{k_0k_1\ell (1-s) n}(k_0k_1\ell (1-s)b n),
\end{align*}
taking expectations in~(\ref{eq_conv_proof}) and
applying~(\ref{bounded_conv}), we obtain~(\ref{eq_conv}) (note also the following elementary fact:
if $\xi\leq b$ a.s., then $\IE\xi\I_A \geq \IE\xi - b\IP[A]$).
\qed

\medskip

Now, we are able to prove Theorem~\ref{t_local_size}.

\medskip
\noindent
{\it Proof of Theorem~\ref{t_local_size}: } 
Consider a closed set $G\subset F^o$ and fix any $\eps>0$.
There exists~$\delta>0$ such that $G\subset F^o_\delta$.
Clearly, for any small enough $\eps'>0$ there exist 
$a_1,\ldots, a_\ell \in(\Q^d\cap F^o_\delta)\setminus\{0\}$ such that
\begin{equation}
\label{covering-}
 \sup_{b\in G} \min_{i=1,\ldots,\ell} \|b-(1-\eps')a_i\| < \eps'/2.
\end{equation}
Let $k_i=\min\{k\in 2\N : k a_i\in 2\Z^d\}$, 
and let $m_i=\max\{m: k_i m \leq (1-\eps')n\}$, $i=1,\ldots,\ell$.
Then, by Lemma~\ref{l_def_beta},
\begin{equation}
\label{snizu}
 \frac{\ln\Eo\eta^0_{k_im_i}(k_ia_im_i)}{k_im_i} \geq \beta(a_i)-\eps
\end{equation}
for all~$n$ large enough, $i=1,\ldots,\ell$.
Using~(\ref{covering-}), we obtain that for any $y\in nG$ there
exists~$i$ such that $\|y-(1-\eps')a_in\| \leq \eps'n/2$.
Then, analogously to the proof of the second claim of Theorem~\ref{t_shape_reach},
we can show that on the event~$M_n^{\delta_0}$
\begin{equation}
\label{ostatok}
 \Eo\eta^{k_ia_im_i}_{n-k_im_i}(y) \geq \delta_0\eps_0^{n-k_im_i}
\end{equation}
for all~$n$ large enough (note that $n-k_im_i\geq \eps'n$).
Now, from~(\ref{supermult}), (\ref{snizu}), and~(\ref{ostatok})
we obtain that 
\begin{align*}
\Eo\eta^0_n(y) &\geq \Eo\eta^0_{k_im_i}(k_ia_im_i)
       \Eo\eta^{k_ia_im_i}_{n-k_im_i}(y)\\
 &\geq \exp\big(k_im_i(\beta(a_i)-\eps)\big)\delta_0\eps_0^{n-k_im_i}\\
 &\geq \delta_0\exp\big(n((1-2\eps')(\beta(a_i)-\eps)-2\eps'\ln\eps_0^{-1})\big).
\end{align*}
  So, from the uniform continuity of~$\beta$ in~$G$ (cf.\ Lemma~\ref{l_convex_beta})
we obtain that
\[
 \liminf_{n\to\infty}
    \min_{x\in nG\cap\Z^d} (n^{-1}\ln \Eo\eta^0_n(x) - \beta(x/n)) \geq 0.
\]

To complete the proof of Theorem~\ref{t_local_size}, we
have to show that 
\begin{equation}
\label{limsup_local_size}
\limsup_{n\to\infty}
    \max_{x\in nG\cap\Z^d} (n^{-1}\ln \Eo\eta^0_n(x) - \beta(x/n)) \leq 0.
\end{equation}

Again, for any small enough $\eps'>0$ there exist 
$a_1,\ldots, a_\ell \in(\Q^d\cap F^o_\delta)\setminus\{0\}$ such that
\begin{equation}
\label{covering+}
 \sup_{b\in G} \min_{i=1,\ldots,\ell} \|b-(1+\eps')a_i\| < \eps'/2.
\end{equation}
Recall $k_i=\min\{k\in 2\N : k a_i\in 2\Z^d\}$, 
and let $m_i'=\min\{m: k_i m \geq (1+\eps')n\}$, $i=1,\ldots,\ell$.

Suppose that there exists $y\in nG$ such that
\[
 \frac{\ln\Eo\eta^0_n(y)}{n} \geq \beta(y/n) + 2\eps.
\]
By~(\ref{covering+}), there
exists~$i$ such that $\|y-(1+\eps')a_in\| \leq \eps'n/2$.
Then, on the one hand, Lemma~\ref{l_def_beta} implies that
\[
 \Eo\eta^0_{k_im_i'}(k_ia_im_i') < \exp\big(k_im_i'(\beta(a_i)+\eps)\big)
\]
for all~$n$ large enough, and, on the other hand,
\begin{align*}
 \Eo\eta^0_{k_im_i'}(k_ia_im_i') &\geq \Eo\eta^0_n(y) 
                         \Eo\eta^y_{k_im_i'-n}(k_ia_im_i'-y)\\
 &\geq \exp\big(n(\beta(y/n)+2\eps)\big)\delta_0\eps_0^{k_im_i'-n}.
\end{align*}
This leads to a contradiction when~$\eps'$ is small enough, 
and thus we obtain~(\ref{limsup_local_size}).
The proof of Theorem~\ref{t_local_size} is completed.
\qed

\medskip
\noindent
{\it Proof of Theorem~\ref{t_rec_trans}: } 
First, note that recurrence implies the existence of ${\mathbf p}$-recurrent
seeds (cf.\ Lemma~3.1 of~\cite{CP} and Definition~\ref{rec_seed} above). 
Such seeds give rise to a supercritical
Galton-Watson process that survives with positive probability (see~\cite{CP} for
details) and so the expected number of the particles at the origin grows exponentially,
thus showing that $\beta(0)>0$. 

On the other hand, if $\beta(0)>0$, then Theorem~\ref{t_local_size}
implies that, with positive $\IP$-probability, there exists $n\geq 1$
(possibly depending on~$\om$) such that $\Eo\eta^0_n(0)>1$. This
implies the existence of a recurrent seed. Indeed,
denote
\[
 \B_\eps(\omega) = \{{\tilde \omega}\in\M : {\tilde \omega}(v)>0 \text{ if and only if }
     \omega(v)>0, |{\tilde \omega}(v)-\omega(v)|<\eps\},
\]
and take $U=\K_{nL_0}$. Choose a small~$\eps$ in such a way that
${\mathtt E}_{{\tilde \omega}} \eta^0_n(0) > 1$ for any~$\tilde\om$ 
such that ${\tilde \omega}_x\in\B_\eps(\omega_x)$ for all $x\in \K_{nL_0}$.
Then, $(\K_{nL_0}, (\B_\eps(\omega_x), x\in \K_{nL_0}))$ is a recurrent seed,
 and so the \BRW{} is recurrent.
\qed

\subsection{Proof of Theorem~\ref{t_complete_shape}}
\label{s_proof_t_complete_shape}
Fix $\eps>0$ and consider any $y\in nG$. By Theorem~\ref{t_local_size},
there exists $n(\om,G)$ (which does not depend on~$y$) such that 
$\Eo\eta^0_n(y)\leq \exp((\beta(y/n)+\eps)n)$ for all $n\geq n(\om,G)$.
We write
\[
 \Po[\eta^0_n(y)\geq \exp((\beta(y/n)+2\eps)n)] \leq
     \frac{\Eo\eta^0_n(y)}{\exp((\beta(y/n)+2\eps)n)}
 \leq e^{-\eps n},
\]
so, by Borel-Cantelli lemma,
\begin{equation}
\label{liminf_complete_shape}
\liminf_{n\to\infty}\max_{x\in nG\cap\Z^d} (n^{-1}\ln \eta^0_n(x) - \beta(x/n)) \leq 0 
    \qquad\text{ $\Po$-a.s.}
\end{equation}

Now, we have to show that 
\begin{equation}
\label{limsup_complete}
\limsup_{n\to\infty}\min_{x\in nG\cap\Z^d} (n^{-1}\ln \eta^0_n(x) - \beta(x/n)) \geq 0 
    \qquad\text{ $\Po$-a.s.}
\end{equation}
In order to prove~(\ref{limsup_complete}), let us first prove that,
for any $a\in\Q^d\cap F^o$
\begin{equation}
\label{liminf_fixed_a}
\liminf_{n\to\infty} \frac{\ln\eta^0_{k_0n}(k_0an)}{k_0n} \geq \beta(a) \qquad\text{ $\Po$-a.s.},
\end{equation}
where $k_0=\min\{k\in 2\N : k a\in 2\Z^d\}$.

\smallskip
\noindent
{\bf Step 1:} First of all, we establish that, for any $\eps>0$
\begin{equation}
\label{Po(liminf)>0}
 \Po\Big[\liminf_{n\to\infty} \frac{\ln\eta^0_{k_0n}(k_0an)}{k_0n} \geq \beta(a)-\eps\Big]>0 
        \qquad\text{ for $\IP$-almost all $\om$}.
\end{equation}
Choose~$\delta>0$ such that $a\in F^o_\delta$, then choose a positive $h\in\Q$ in such a way that
$a(1-h)^{-1}\in F^o_\delta$, then let $k_1=\min\{k\in 2\N : k a\in 2\Z^d, kh\in 2\Z\}$.
Abbreviate
\[
 g_n = \IP[k_1 a n \in R_\omega^{0,\delta}(k_1(1-h)n)];
\]
recall that, by~(\ref{incl_Rn}), $g_n\to 1$ as $n\to\infty$.
By virtue of~(\ref{bounded_conv}), 
\[
 \frac{\IE\big(\1{k_1 a n \in R_\omega^{0,\delta}(k_1(1-h)n)}\ln\Eo\eta^0_{k_1(1-h)n}(k_1an)\big)}
         {k_1(1-h)n} \to \beta((1-h)^{-1}a),
\]
so one can choose~$n_1$ such that 
\begin{equation}
\label{choosing_n_1}
\frac{\IE\big(\1{k_1 a n_1 \in R_\omega^{0,\delta}(k_1(1-h)n_1)}\ln\Eo\eta^0_{k_1(1-h)n_1}(k_1an_1)\big)}
         {k_1(1-h)n_1} \geq  \beta((1-h)^{-1}a) - \eps,
\end{equation}
and also
\begin{equation}
\label{1st_restr_n_1}
 1-2h < (1-h)g_{n_1}+\|a\|(1-g_{n_1}) < 1,
\end{equation}
\begin{align}
 (\beta((1-h)^{-1}a)-2\eps)(1-h) & < (\beta((1-h)^{-1}a)-\eps)(1-h)g_{n_1} \nonumber\\
 & \qquad {} - \|a\|\ln\eps_0^{-1}(1-g_{n_1}). \label{2nd_restr_n_1} 
\end{align}

Now, we construct a branching process in random environment $(\Upsilon_\ell,\ell=0,1,2,\ldots)$
in the following way. Here, $\Upsilon_\ell$ stands for the size of $\ell$th generation of this 
process. With respect to the original process, the particles of $\ell$th generation are
in $k_1 a \ell n_1$ at time~$t_\ell$ defined below (note that it means that
$\Upsilon_\ell \leq \eta^0_{t_\ell}(k_1 a \ell n_1)$, but the equality should not
necessarily hold true, there may be also some particles in $k_1 a \ell n_1$ at time~$t_\ell$
which do not belong to this branching process in random environment). Specifically,
the initial particle is considered the particle of $0$th generation, and we set 
$\Upsilon_0=1$, $t_0=0$. Inductively, consider the $\Upsilon_{\ell-1}$ particles
of $(\ell-1)$th generation, situated in $k_1 a (\ell-1) n_1$ at the moment~$t_{\ell-1}$.
Then, the particles of $\ell$th generation are their descendants 
% of the particles of $(\ell-1)$th generation 
which are in $k_1 a \ell n_1$ at time~$t_\ell$, where
\[
 t_\ell = \left\{
 \begin{array}{ll}
 t_{\ell-1}+k_1(1-h)n_1, & \text{on }\{k_1 a \ell n_1\in 
              R_\omega^{k_1 a (\ell-1) n_1,\delta}(k_1(1-h)n_1)\}\\
 t_{\ell-1}+k_1\|a\|n_1, & \text{on }\{k_1 a \ell n_1\notin 
              R_\omega^{k_1 a (\ell-1) n_1,\delta}(k_1(1-h)n_1)\}.
\end{array}
\right.
\]
By Birkhoff's theorem,
\[
 \lim_{m\to\infty} \frac{t_m}{m} = k_1n_1 \big((1-h)g_{n_1}+\|a\|(1-g_{n_1})\big),
\]
and, by Condition~UE and~(\ref{2nd_restr_n_1}),
\begin{align}
\IE \ln \Eo\Upsilon_1 &\geq (\beta((1-h)^{-1}a)-\eps)k_1(1-h)n_1g_{n_1} -
      k_1\|a\|n_1\ln\eps_0^{-1}(1-g_{n_1}) \nonumber\\
 &\geq (\beta((1-h)^{-1}a)-2\eps)k_1(1-h)n_1.
\label{expect_bpre}
\end{align}
Assume that~$\eps$ is so small that $\beta((1-h)^{-1}a)>2\eps$.
Since (by Condition~UE and using the fact that $a(1-h)^{-1}\in F^o_\delta$)
\[
 \Po[\Upsilon_1\geq 1] \geq \min\{\eps_0^{k_1\|a\|n_1},\delta^{k_1(1-h)n_1}\},
\]
one can use e.g.\
Theorem~1 of~\cite{AK}, or Theorem~5.5 and Proposition~6.2 of~\cite{T}
to obtain that
\begin{equation}
\label{surv_bpre}
 \Po[\text{the process $\Upsilon$ survives}]>0 \qquad \text{ $\IP$-a.s.},
\end{equation}
and
\begin{equation}
\label{growth_bpre}
 \Po[\liminf_{m\to\infty} m^{-1}\ln\Upsilon_m\geq (\beta((1\!-\!h)^{-1}a)-2\eps)k_1(1\!-\!h)n_1
   \mid \Upsilon\text{ survives}] = 1.
\end{equation}

On the event $\{\liminf_{m\to\infty} m^{-1}\ln\Upsilon_m\geq (\beta((1-h)^{-1}a)-2\eps)k_1(1-h)n_1\}$
one can choose~$m_0$ (depending on~$\om$) such that 
\[
 \frac{\ln\Upsilon_m}{m} \geq (\beta((1-h)^{-1}a)-3\eps)k_1(1-h)n_1
\]
and
\[
 (1-2h)k_1n_1 \leq \frac{t_m}{m} \leq k_1n_1
\]
for all $m\geq m_0$. Then, at the moment~$t_m$ we have at least
\[
\Upsilon_m\geq \exp\big(mk_1n_1(1-h)(\beta((1-h)^{-1}a)-3\eps)\big)
\]
particles in $k_1an_1m$. By Condition~UE, each of those particles has a descendant in
$k_1an_1m$ at time $k_1n_1m$ with probability at least $\eps_0^{2hk_1n_1m}$.
So, using the large deviation bound for the binomial distribution (cf.\ e.g.\
formula~(34) of~\cite{CP}), we obtain that for some positive~$C_1,C_2$
\begin{align*}
 \lefteqn{\Po[\eta^0_{k_1n_1m}(k_1an_1m) \geq \exp\big(mk_1n_1((1-h)(\beta((1-h)^{-1}a)-3\eps)
      -2h\ln\eps_0^{-1})\big)]}~~~~~~~~~~~~~~~~~~~~~~~~~~~~~~~~~~~~~~~~~~~~~~~~~~~~~~~~~~~~~~~~~~~~~~ \\
 &\geq 1-\exp(-C_1e^{C_2m}).
\end{align*}
 Using Condition~UE again, we obtain~(\ref{Po(liminf)>0}).

\smallskip
\noindent
{\bf Step 2:}
Now, let us show that~(\ref{Po(liminf)>0}) implies~(\ref{liminf_fixed_a}).
This is easy in the case when the \BRW{} is recurrent. Indeed, in this case
it can be shown that a.s.\ the origin will be visited by infinitely many ``free''
particles (i.e., none of them is a descendant of another; see~\cite{CP} for more details).
Each of those particles gives rise to a copy of the 
branching process in random environment constructed above
(they use the same environment, but are conditionally independent); so, with probability~$1$
at least one of them survives, and from this we obtain~(\ref{liminf_fixed_a})
in the recurrent case. However, this argument does not work in the case when the \BRW{} is
transient, so we present a general argument that works in both cases.

Abbreviate ${\hat u}=\max\{1-h,\|a\|\}$. Let 
\[
 S_x = \{y\in\Z^d : \text{there is }\ell\in\Z_+ \text{ such that }
                  \|x+\ell k_1an_1-y\|\leq L_0k_1{\hat u}n_1+\rho\}.
\]
The key observation is that the branching process~$\Upsilon$ constructed above
depends only on the environment inside $S_0$. In particular the probability in 
(\ref{surv_bpre}) only depends on $\omega_x, x \in S_0$. 
Suppose that~$\Upsilon_m=0$ for
some~$m$; this means that up to (``real'') time $k_1{\hat u}n_1m$
the branching process in random environment became extinct.

Then, the idea is
the following: with positive probability 
it happens that a particle goes outside the already explored
part of~$\om$ ``without revealing more environment'',
 and we can construct a new branching process
in random environment, independent of the previous one. 

At the moment $k_1{\hat u}n_1m$ (when we know that the ``initial'' 
branching process~$\Upsilon$ in random environment became extinct),
let us remove all particles except one from the process~$\eta$ (for definiteness,
choose the remaining particle uniformly among the particles that are present
at time $k_1{\hat u}n_1m$). 
Let~$z$ be the ``initial'' (i.e., at time $k_1{\hat u}n_1m$) position of this particle,
note that $\|z\|\leq L_0k_1{\hat u}n_1m$.
We let this particle perform the uniform induced random walk~$\xi^z$
(i.e., immediately removing from the process other particles that may eventually appear)
until some random moment~$\tau$ defined below.
Let
\begin{align*}
 \Gamma(\ell) &=\big\{x\in\Z^d : \text{there exists }y\in\Z^d \text{ such that }y\in S_x\text{ and}\\
& \qquad \qquad \text{ either }\|y\|\leq L_0 k_1{\hat u}n_1m
 \text{ or } y\in\{\xi^z_0,\ldots,\xi^z_\ell\}\big\}
\end{align*}
be the set of sites from where the construction analogous to the 
construction of the above branching process~$\Upsilon$
may depend on already revealed pieces of the environment.

\begin{figure}
\centering
\includegraphics{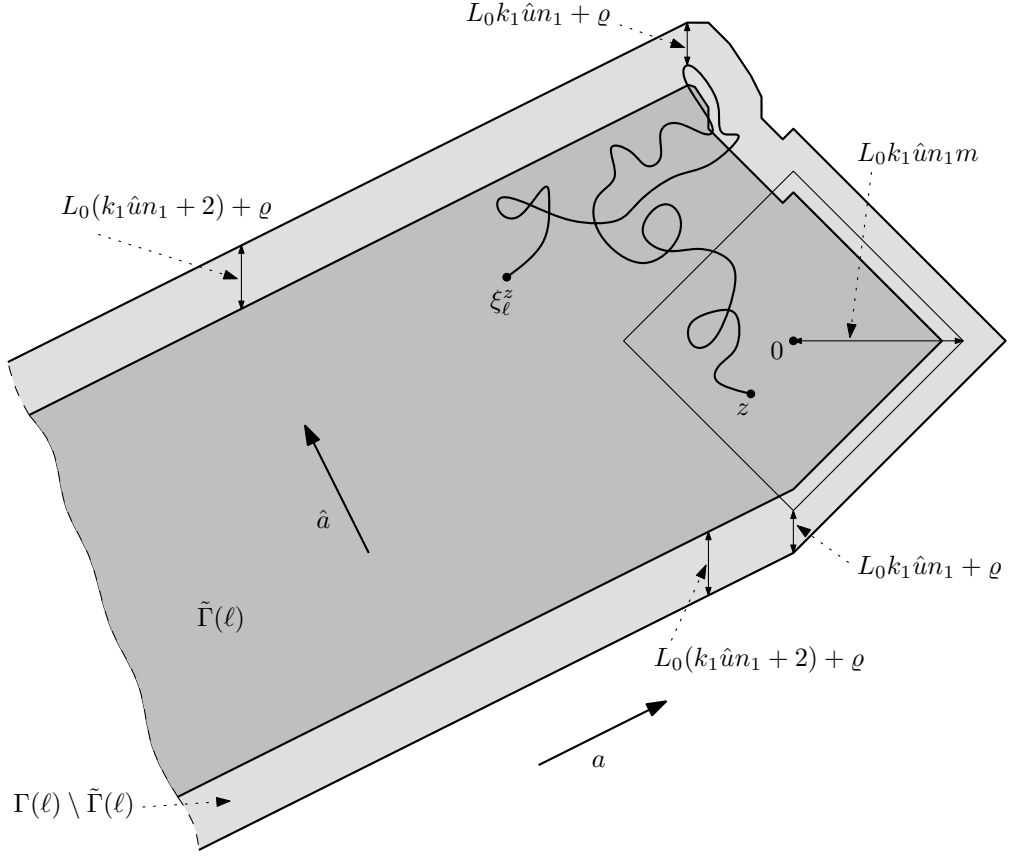}
\caption{On the definition of the sets $\Gamma(\ell)$ and ${\tilde\Gamma}(\ell)$
(recall that $\|\cdot\|$ is the ${\mathcal L}_1$-norm)}
\label{df_Gamma}
\end{figure}

Recall the notation $Z_t$ from Section \ref{s_ind_rw}. We
define
\begin{align*}
 \tau &=\min\{s\in\N : \text{ there exists }s'<s \\
 &\qquad \qquad \text{ such that } Z_{s'+1}=\ldots=Z_s=1
   \text{ and } \xi^z_s\notin \Gamma(s')\},
\end{align*}
and also (see Figure~\ref{df_Gamma})
\[
 {\tilde\Gamma}(\ell)=\big\{y\in\Gamma(\ell) : 
            \dist(\{y\},\Z^d\setminus \Gamma(\ell))>L_0(k_1{\hat u}n_1+2)+\rho\big\}.
\]
Clearly, if $\xi^z_\ell \in {\tilde\Gamma}(\ell)$, then $\Gamma(\ell+1)=\Gamma(\ell)$.
 From Condition~UE it follows that a.s.\ there exists an 
infinite sequence $(\tau'_k,k=1,2,3,\ldots)$ such that 
$\xi^z_{\tau'_k}\in \Gamma(\tau'_k) \setminus {\tilde\Gamma}(\tau'_k)$ for all~$k$. 
Indeed, if ${\hat a}\in\R^d\setminus\{0\}$ is such that ${\hat a}\cdot a = 0$,
then the fact that $|(\xi^z_{t_1}-\xi^z_{t_2})\cdot {\hat a}|$ is large enough
guarantees that there is~$k$ such that $\tau'_k\in [t_1,t_2]$.
Since from any $u\in \Gamma(\ell) \setminus {\tilde\Gamma}(\ell)$ the particle
can perform any given $L_0(k_1{\hat u}n_1+2)+\rho$ unit steps with $Z$-value $1$ with
uniformly positive probability (and so it can go out of $\Gamma(\ell)$ with
at least that probability without revealing more environment), this shows
that $\tau<\infty$ a.s. 

Now, at time $k_1{\hat u}n_1m + \tau$ we can start 
an independent copy~$\Upsilon'$ of the branching process in random environment~$\Upsilon$
constructed above. If it happens that~$\Upsilon'$ dies out as well, repeating this 
construction one can start another independent copy~$\Upsilon''$, and so on.
Eventually, one of those branching processes in random environment will
survive, and so~(\ref{liminf_fixed_a}) follows from~(\ref{growth_bpre}).

\smallskip
Now, let us prove~(\ref{limsup_complete}). We use the same method as in the proof of
the first part of Theorem~\ref{t_local_size}.
Namely, suppose that $\delta>0$ is such that $G\subset F^o_\delta$.
For any small enough $\eps'>0$ there exist 
$a_1,\ldots, a_\ell \in(\Q^d\cap F^o_\delta)\setminus\{0\}$ 
such that~(\ref{covering-}) holds.
Let $k_i=\min\{k\in 2\N : k a_i\in 2\Z^d\}$, 
and let $m_i=\max\{m: k_i m \leq (1-\eps'n)\}$, $i=1,\ldots,\ell$.
By~(\ref{liminf_fixed_a}), provided that $\eps,n_1,h$ are chosen in such a way 
that 
\[
 (1-h)(\beta((1-h)^{-1}a_i)-3\eps)-2h\ln\eps_0^{-1} \geq \beta(a_i)-\eps''
\]
we have
\begin{equation}
\label{snizu_complete}
 \frac{\ln\eta^0_{k_im_i}(k_ia_im_i)}{k_im_i} \geq \beta(a_i)-\eps''
\end{equation}
for all~$n$ large enough, $i=1,\ldots,\ell$.
Using~(\ref{covering-}), we obtain that for any $y\in nG$ there
exists~$i$ such that $\|y-(1-\eps')a_in\| \leq \eps'n/2$.
Then, we can show that on the event~$M_n^{\delta_0}$
\begin{equation}
\label{prob_omega_reach}
 \Po[\eta^{k_ia_im_i}_{n-k_im_i}(y) \geq 1] \geq \delta_0\eps_0^{n-k_im_i}.
\end{equation}
So, since the particles that are in $k_ia_im_i$ at time $k_im_i$
act independently, again using a large deviation bound for the binomial distribution
together with~(\ref{snizu_complete}) and~(\ref{prob_omega_reach}),
we obtain that for some positive constants $C_3,C_4$
\[
 \Po[\eta^0_n(y)\geq \exp\big(k_im_i(\beta(a_i)-2\eps)\big)] \geq 1 - \exp(-C_3e^{C_4n}),
\]
and this proves~(\ref{limsup_complete}). Theorem~\ref{t_complete_shape}
now follows from~(\ref{liminf_complete_shape}) and~(\ref{limsup_complete}).
\qed

\subsection{Exponential growth of the total number of particles and 
            proof of Theorem~\ref{t_total_size}.} 
\label{s_proof_exp_growth}
First, we prove the property~(\ref{sup_beta>0}).
For the case of recurrent \BRW, (\ref{sup_beta>0}) immediately follows from Theorem~\ref{t_rec_trans},
so we only need to consider the case of transient \BRW. Without restriction of
generality, one can assume that any particle can generate at most two particles
at the next moment, and that there is~$h>0$ such that the probability
of generating two offspring is (depending on the location) $0$ or~$h$, i.e.,
\begin{equation}
\label{reduce_brw}
\IP[\omega_0(|v|\leq 2)=1]=1, \qquad \IP[\omega_0(|v|=2)\in\{0,h\}]=1.
\end{equation}
This is no restriction 
since it can be easily shown that any \BRW{} in random environment
(satisfying Condition~B)
dominates some \BRW{} in random environment satisfying~(\ref{reduce_brw}).
Also, it is clear that the latter \BRW{} in random environment can be defined
in such a way that Condition~UE still holds (possibly with another constant;
however, to keep the notations simple, we will assume in our argument that
Condition~UE holds with $\eps_0$).

In the situation~(\ref{reduce_brw}), if $\omega_x(|v|=2)=h$, we say that~$x$ is a branching site,
in the case when $\omega_x(|v|=2)=0$, we say that~$x$ is a non-branching site.

Then, we can  make a further simplifying assumption: we suppose
that \emph{immediate} descendants of a particle jump independently
(so that we are in the situation considered e.g.\ in~\cite{CMP,MP1,MP2}),
which means that, for a particle in~$x$, we  first decide if the particles
generate 1 or 2 offspring, and then each of these offspring jump 
independently with probabilities $(\Po[\xi^x_1=y], y\in x+\A)$. 
Here, $\xi^x$ is the uniform induced random walk starting from~$x$.
% ;
% note also that, when~(\ref{reduce_brw}) holds, the uniform induced random walk
% and the particle-uniform induced random walk have the same law.
The reason why we can make this assumption without loss of generality
is that the two \BRW s are in some sense equivalent, as mentioned  in Remark 
\ref{rem_complete_shape}. 

Denote by 
\[
 \psi_\omega(x) = \Po[\eta^x_n(x)=0 \text{ for all }n\geq 1]
\]
the probability that none of the descendants of a particle in~$x$ ever comes back.
In the next lemma we prove that $\psi_\omega(x)$ is uniformly positive.

\begin{lmm}
\label{l_noreturn}
There exists $\theta>0$, depending only on~$h$ and~$\eps_0$, such that,
for any transient \BRW{} satisfying~(\ref{reduce_brw}) and with independent
immediate descendants, we have $\psi_\omega(x)>\theta$, for all~$x\in\Z^d$
and for $\IP$-almost all~$\om$.
\end{lmm}

\noindent
{\it Proof of Lemma \ref{l_noreturn}: } 
By stationarity, it is enough to prove that $\IP[\psi_\omega(0)>\theta]=1$.
We argue by contradiction.
The idea is  that, when $\psi_\omega(0)$ is small enough, 
 the \BRW{} should be recurrent. We split the space of environments in three 
parts, and 
consider three cases accordingly.

\medskip
\noindent
\underline{Case 1.} Suppose that $\omega_0(|v|=2)=h$, i.e., $0$ is 
a branching site. Abbreviating $b:=\Po[\eta^0_n(0)=0 \text{ for all }n\geq 1\mid \ZZ^0_1=1]$,
we have $\psi_\omega(0) = (1-h)b+hb^2$, so 
\[
 \Po[\eta^0_n(0)=0 \text{ for all }n\!\geq \! 1\mid \ZZ^0_1\!=\!1] = 
     \frac{-(1\!-\!h)+\sqrt{(1\!-\!h)^2+4h \psi_\omega(0)}}{2h}.
\]
Thus, $\Po[\eta^0_n(0)=0 \text{ for all }n\geq 1\mid \ZZ^0_1=1]$ can be made
arbitrarily close to~$0$ by making $\psi_\omega(0)$ small. Now, if
\[
 \Po[\eta^0_n(0)=0 \text{ for all }n\geq 1\mid \ZZ^0_1=1] < 
\frac{h}{1+h},
\]
then, as in the proof of Theorem~\ref{t_rec_trans}, it is straightforward
to show that there exists a recurrent seed and so
the \BRW{} is recurrent.

\medskip
\noindent
\underline{Case 2.} We suppose that $\omega_0(|v|=2)=0$, but there is a branching site
in $\K_{\rho-1}\setminus\{0\}$. In this case we can choose a branching site there which
is closest to the origin, so that it 
is accessible from the origin by a path of non-branching sites. This means that there exist
$k\leq \rho-1$, $0=x_0,x_1,\ldots,x_k$ such that $\|x_{i+1}-x_i\|=1$,
$i=0,\ldots,k$, $x_1,\ldots,x_{k-1}$ are non-branching, and~$x_k$ is a
branching site (and we have also $\|x_k\|=k$).

Then, again we obtain that if $\psi_\omega(0)$ is too small, then the \BRW{}
should be recurrent. For this, proceed as follows. First, define
\[
 \tau_1 = \min\{n\geq 1 : \eta^0_n(0)\geq 1\}\;,
\]
so that $\Po[\tau_1<\infty]=1-\psi_\omega(0)$.
At the moment~$\tau_1$ (provided $\tau_1<\infty$), consider one of the particles which
are in~$0$ and let
\begin{align*}
 \tau_2 &=\min\{n\geq\tau_1+1 : \text{ at least one of the descendants}\\
   &\qquad \qquad \qquad \qquad \qquad \text{ of that particle is in~$0$
         at time~$n$}\},
\end{align*}
then, repeat this procedure to define $\tau_3$ on $\{\tau_2<\infty\}$,
$\tau_4$ on $\{\tau_3<\infty\}$, and so on. Clearly, $\Po[\tau_m<\infty]=(1-\psi_\omega(0))^m$.

Being $\xi^0$ the uniform induced random walk starting from~$0$, define the event
\[
 A = \{\xi^0_i=x_i, i=1,\ldots,k, \xi^0_i=x_{2k-i}, i=k+1,\ldots,2k, \ZZ^0_{k+1}=2\},
\]
i.e., the initial particle goes straight to~$x_k$, branches there, and then the ``first''
descendant goes straight back to the origin (note that we do not assume anything about the
second descendant). Clearly, we have
\[
 \Po[A\mid \tau_1<\infty] \geq \frac{h\eps_0^{2k}}{1-\psi_\omega(0)}
                             \geq h\eps_0^{2(\rho-1)}.
\]

On the event~$A$, the second particle (generated in $x_k$) goes to~$0$
with probability at least $\eps_0^{\rho-1}$. Take~$m$ such that
$mh\eps_0^{3(\rho-1)}>1$: if $\psi_\omega(0)$ is small enough,
then $(1-\psi_\omega(0))^m mh\eps_0^{3(\rho-1)}>1$. Then we obtain
that there exists large enough~$T$ (which actually depends on~$\om$)
such that up to time~$T$ the mean number of ``free'' particles
 that visit the origin
is greater than~$1$. (Recall that, by ``{\em free}'' particles we mean that none of them is the descendant of
another one, see Definition~\ref{rec_seed}.)
As above,
this implies the existence of recurrent seeds and so
the \BRW{} is recurrent.

\medskip
\noindent
\underline{Case 3.} Suppose that $\omega_x(|v|=2)=0$ for all $x\in\K_{\rho-1}$,
so that there are no branching sites in $\K_{\rho-1}$.
Let
\[
 \partial_e\K_{\rho-1} = \big\{x\in \K_{\rho-1}^c: \dist(\{x\},\K_{\rho-1})
\leq L_0\big\}
\]
be the annulus (or extended external boundary) of $\K_{\rho-1}$.
Denote by
\[
 g_\omega(x,B) = \Po\Big[\text{there exists } n\geq 0\text{ such that }
                     \sum_{y\in B}\eta^x_n(y)\geq 1\Big]
\]
the probability of ever having particles in~$B\subset\Z^d$
for the process starting from~$x$. For any $x\in\partial_e\K_{\rho-1}$,
we see from Condition~UE that 
$\psi_\omega(0) > \eps_0^{\rho+L_0}(1-h)^{L_0}(1-g_\omega(x,\partial_e\K_{\rho-1}))$,
so
\begin{equation}
\label{inf_g(xK)}
 \inf_{x\in\partial_e\K_{\rho-1}} g_\omega(x,\partial_e\K_{\rho-1}) 
             > 1-\eps_0^{-(\rho+L_0)}(1-h)^{-L_0}\psi_\omega(0).
\end{equation}
At this point, the important observation is that, 
for any $x\in\partial_e\K_{\rho-1}$, the quantity
$g_\omega(x,\partial_e\K_{\rho-1})$ is a function of 
the environment outside $\K_{\rho-1}$,
and so it is independent of $\omega_0$ by Condition~I. 
So, supposing that $\psi_\omega(0)$ is small enough,
one can construct a recurrent seed as follows.
Define
\[
 g^{(k)}_\omega(x,B) = \Po\Big[\text{there exists } 0\leq n\leq k\text{ such that }
                     \sum_{y\in B}\eta^x_n(y)\geq 1\Big],
\]
and note that $g^{(k)}_\omega(x,B)\to g_\omega(x,B)$ as $k\to\infty$.
Then, suppose that $\psi_\omega(0) < {\tilde\delta}$, where (a small number) ${\tilde\delta}$
is such that 
\begin{equation}
\label{eq_df_tilde_delta}
 \frac{\eps_0^{\rho-1}}{\eps_0^{\rho-1}+\eps_0^{-(\rho+L_0)}(1-h)^{-L_0}{\tilde\delta}} 
             > (1+h)^{-1}.
\end{equation}
Choose first large enough~$k_0$ and then small enough~$\eps$ in such a way that
\begin{equation}
\label{choose_k0}
 \inf_{x\in\partial_e\K_{\rho-1}} g^{(k_0)}_{\tilde\omega}(x,\partial_e\K_{\rho-1}) 
             > 1-\eps_0^{-(\rho+L_0)}(1-h)^{-L_0}{\tilde\delta}
\end{equation}
for all $\tilde\om$ such that ${\tilde\omega}_x\in\B_\eps(\omega_x)$ 
 for all $x\in\K_{\rho+k_0L_0}\setminus\K_{\rho-1}$
(here, $\B_\eps(\cdot)$ denotes the ball of radius~$\eps$ in~$\M$ 
with respect to any fixed metrics). By the contradiction assumption and 
by (\ref{inf_g(xK)}), the set of such  $\tilde\om$'s has positive 
$\IP$-probability.

Then, the idea is the following:
put a branching site at the origin and suppose that in a sufficiently
large region around~$0$ (excluding $\K_{\rho-1}$) the environment is ``close''
to the environment above. More precisely, we consider the $(U,\HH)$-seed
with $U=\K_{\rho+k_0L_0}$ and
\[
 \HH_x = \left\{\begin{array}{ll}
           \{\omega\in\M : \omega(|v|=2)=h\}, & \text{for }x=0,\\
           \M, & \text{for }x\in \K_{\rho-1}\setminus\{0\},\\
           \B_\eps(\omega_x), & \text{for }x\in \K_{\rho+k_0L_0}\setminus\K_{\rho-1}.
        \end{array} \right.
\]
Condition~I implies that $\IP[\tilde\om : {\tilde\omega}_x\in\HH_x \text{ for all }x\in U]>0$.

 Now, from each site in~$\K_{\rho-1}$, the uniform induced random walk goes
to~$0$ without leaving~$\K_{\rho-1}$ with probability at least $\eps_0^{\rho-1}$,
and any particle in $\partial_e\K_{\rho-1}$ sends at least one descendant 
to 
%%$\K_{\rho-1}$
%% F:
$\partial_e\K_{\rho-1}$
 at least with probability given by~(\ref{choose_k0}).
 Suppose without restriction of generality that $\rho-1 \geq L_0$.
Then,
the probability that any particle starting from~$\A$ sends at least one descendant to~$0$
is at least\footnote{Compare with the following situation. There are two coins, for the
coin 1 the probability of head is~$p$, for the coin 2 the probability of head is~$q$.
We flip the coins in an alternate fashion, i.e., 1,2,1,2,1,2,\dots Then, the probability
that the first head comes from the coin 1 is equal to $\frac{p}{p+q-pq}>\frac{p}{p+q}$.}
\[
 \frac{\eps_0^{\rho-1}}{\eps_0^{\rho-1}+\eps_0^{-(\rho+L_0)}(1-h)^{-L_0}{\tilde\delta}},
\]
which is greater than $(1+h)^{-1}$ by~(\ref{eq_df_tilde_delta}). So, we obtain 
a recurrent seed.

This concludes the proof of Lemma~\ref{l_noreturn}.
\qed

We continue the proof of Theorem \ref{t_total_size}.
Let $U_n=\{\xi^0_0,\xi^0_1,\ldots,\xi^0_n\}$ be the range of the uniform
induced random walk up to time~$n$. We now prove that, due to Lemma~\ref{l_noreturn},
$|U_n|$ is of order~$n$ with at least constant probability. Let 
\[
N(x,y)=\sum_{i=1}^\infty \1{\xi^x_i=y}
\]
be the number of times that the uniform induced 
random walk starting from~$x$ visits~$y$.
 Using Lemma~\ref{l_noreturn}, we obtain that $\Eo N(x,x) \leq \theta^{-1}$
for all $x\in\Z^d$ and for $\IP$-almost all~$\om$. 
On the other hand
\[
 n\leq \sum_{x\in U_n} N(0,x),
\]
so, taking expectation,
\begin{align*}
 n &\leq \sum_{x\in\Z^d}\Eo(\1{x\in U_n} N(0,x)) \\
 & = \sum_{x\in\Z^d} \Eo(\1{x\in U_n} N(0,x)\mid x\in U_n)\Po[x\in U_n]\\
 &\leq (1+\theta^{-1})\Eo |U_n|.
\end{align*}
Since trivially $|U_n|\leq n$, using the fact that
for any random variable~$X$ with $0\leq X\leq a$ a.s.\ and
$\Eo X\geq b$ it is true that $\Po[X\geq b/2]\geq b/(2a)$, we obtain that,
for $\IP$-almost all~$\om$,
\begin{equation}
\label{range>Cn}
 \Po[|U_n|\geq (1+\theta^{-1})^{-1}n/2] \geq (1+\theta^{-1})^{-1}/2.
\end{equation}

\begin{figure}
\centering
\includegraphics{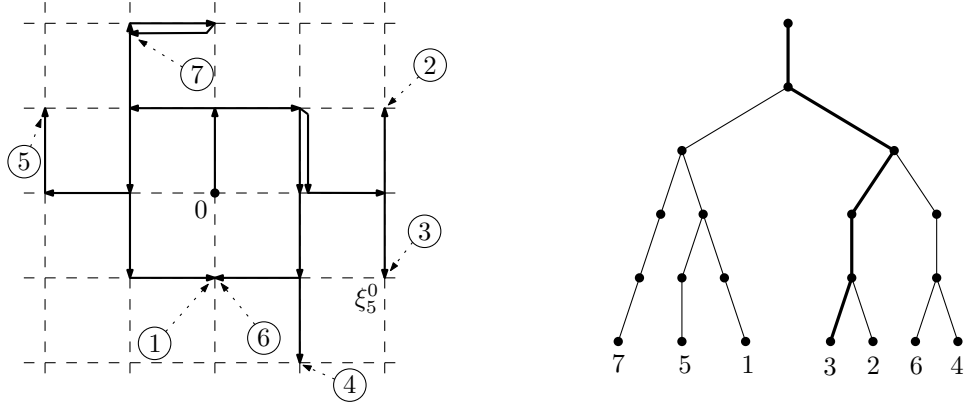}
\caption{A realisation of the \BRW{} 
(together with the uniform induced random walk) up to time $n=5$,
and the corresponding genealogical tree. 
We have $\zeta_5=3$, $\br(i,5)=3$ for $i=1,\ldots,6$, $\br(7,5)=2$.}
\label{f_gen_tree}
\end{figure}

Consider the evolution of the \BRW{} up to time~$n$, and let us enumerate
the~$\ZZ^0_n$ particles that are present at time~$n$ in random order
(i.e., select one particle at random and attach the label ``1'' to it, then 
select one of the unlabelled particles and put the label ``2'' to it, and so on).
We define $\br(i,n)$ to be the number of bifurcations on the path from the root 
to the particle (of the $n$th generation) labelled~$i$ 
on the genealogical tree of the \BRW{} (see Figure~\ref{f_gen_tree}).
Let~$\zeta_n$ be the label assigned to the particle corresponding to the uniform
induced random walk at time~$n$; clearly, given the realisation of the 
genealogical tree,
 \begin{equation}
\label{eq_prob_zeta}
\{\zeta_n=i\} \; \mbox{ has probability } 2^{-\br(i,n)}.
\end{equation}

Let us prove now that, on the event that the range of the uniform induced random walk
is linear in~$n$, with large probability $\br(\zeta_n,n)$ will be linear as well.
To this end, we construct a set $\Phi\subset\Z_+$ in the following way:
 first, we have $0\in\Phi$.
Inductively, suppose that the set $\Phi\cap\{0,\ldots,k-1\}$ was already constructed.
Then, $k\in\Phi$ if and only if the following holds:
\begin{itemize}
 \item there exists $y\in\Z^d$ such that $\|\xi^0_k-y\|=\rho$, and
 \item $\|\xi^0_m-y\|\geq \rho$ for all $m\leq k$.
%, and
% \item $k-k'\geq \rho+1$ for all $k'\in \Phi\cap\{0,\ldots,k-1\}$.
\end{itemize}
Define the cubes
\[
 \hK_m = \{x\in\Z^d : \max(|x^{(1)}|,\ldots,|x^{(d)}|)\leq m\},
\]
note that
\[
 \Z^d = \bigcup_{z\in\Z^d} ((2\rho+2L_0+1)z+\hK_{\rho+L_0}).
\]
Let $t_0=0$, $z_0=0$. Inductively, define for $i\geq 1$
\[
 t_i = \min\Big\{ t : \xi^0_t\notin \bigcup_{j=0}^{i-1} ((2\rho+2L_0+1)z_j+\hK_{\rho+L_0})\Big\},
\]
and $z_i$ is such that $\xi^0_{t_i}\in (2\rho+2L_0+1)z_i+\hK_{\rho+L_0}$.
Observe that, for any~$y\in\Z^d$ such that 
$\dist(\{y\},\Z^d\setminus ((2\rho+2L_0+1)z_i+\hK_{\rho+L_0}))\leq L_0$
we have
\[
 \big(y+(\K_\rho\setminus\K_{\rho-1})\big) \setminus 
  \Big(\bigcup_{x\notin(2\rho+2L_0+1)z_i+\hK_{\rho+L_0}}(x+\K_{\rho-1})\Big) \neq \emptyset, 
\]
so $t_i\in\Phi$ for all $i\geq 0$. 
Since $|\hK_{\rho+L_0}|=(2\rho+2L_0+1)^d$, we have
\begin{equation}
\label{Phi_range}
 |\Phi\cap\{0,\ldots,n\}| \geq (2\rho+2L_0+1)^{-d}|U_n|.
\end{equation}

Let ${\mathcal F}_m$ be the sigma-algebra generated by 
$(\xi^0_k,\br(\zeta_k,k),k\leq m)$.
By definition of the set~$\Phi$,
\begin{equation}
\label{annealed_next_branching}
 \PA[\br(\zeta_{m+\rho},m+\rho)\geq \br(\zeta_m,m)+1 \mid {\mathcal F}_m]
      \geq h\eps_0^\rho \IP[\omega_0(|v|=2)=h].
\end{equation}
Write $\Phi=\{0=\phi_0,\phi_1,\phi_2,\ldots\}$, where $\phi_{i+1}\geq \phi_i$ for all $i\geq 0$.
Abbreviate $\alpha_1=\frac{(1+\theta^{-1})^{-1}}{2}$, $\alpha_2=(2\rho+2L_0+1)^{-d}$,
$\alpha_3=h\eps_0^\rho \IP[\omega_0(|v|=2)=h]$, and $\alpha_4 = \frac{1}{2}\alpha_1\alpha_2\alpha_3\rho^{-1}$.
The property~(\ref{annealed_next_branching}) implies that, for some $C_1>0$,
\[
 \PA[\br(\zeta_{\phi_{\lf\alpha_1\alpha_2n\rf}},\phi_{\lf\alpha_1\alpha_2n\rf})\geq \alpha_4 n] \geq 
   1-e^{-C_1n}.
\]
So, using (\ref{range>Cn}) and (\ref{Phi_range})
together with the elementary inequality $\PA[A\mid B] \geq 1-\frac{\PA[A^c]}{\PA[B]}$,
we obtain for some $C_2>0$ that
\begin{align*}
\PA[\br(\zeta_n,n)\geq\alpha_4 n] &\geq \PA[\br(\zeta_n,n)\geq\alpha_4 n\mid |U_n|\geq \alpha_1n]
        \PA[|U_n|\geq \alpha_1n]\\
  &\geq \PA[\br(\zeta_{\phi_{\lf\alpha_1\alpha_2n\rf}},\phi_{\lf\alpha_1\alpha_2n\rf})\geq \alpha_4 n
        \mid |U_n|\geq \alpha_1n]    \\
           & \quad {}\times    \PA[|U_n|\geq \alpha_1n]\\
 &\geq C_2.
\end{align*}
Therefore, with $\IP$-probability at least $C_3=C_2/2$, 
\[
 \Po[\br(\zeta_n,n)\geq\alpha_4 n] \geq C_3 \, .
\]

Let ${\hat{\mathcal F}}_n$ be the sigma-algebra generated by the \BRW{}
up to time~$n$. 
Applying~(\ref{eq_prob_zeta}), we obtain that, for such $\om$'s,
\[
 \Po[\br(\zeta_n,n)\geq\alpha_4 n\mid {\hat{\mathcal F}}_n]
= \sum_{i:\br(i,n)\geq\alpha_4 n} 2^{-\br(i,n)}
\leq 2^{-\alpha_4n}\ZZ^0_n\,,
\]
so, taking expectations, $\Eo \ZZ^0_n \geq C_3 2^{\alpha_4n}$.
As observed in the beginning of this section, this proves~(\ref{sup_beta>0}).

Now, it remains only to prove~(\ref{eq_total_size}).
Since
\[
 \ln \max_{y\in \K_{L_0n}}\eta^0_n(y) \leq \ln\ZZ^0_n \leq d\ln(2L_0n+1) 
          + \ln \max_{y\in \K_{L_0n}}\eta^0_n(y)
\]
(note that $|\K_{L_0n}|<(2L_0n+1)^d$),
the property~(\ref{eq_total_size}) follows from Theorem~\ref{t_complete_shape}.
\qed

\end{document}